\documentclass[11pt
]{amsart}
\usepackage{
amssymb
}
\newcommand{\no}[1]{#1}
\renewcommand{\no}[1]{}

\no{\usepackage{times}\usepackage[
subscriptcorrection, slantedGreek, 
nofontinfo]{mtpro}
\renewcommand{\Delta}{\upDelta}
}
\date{\today}

\numberwithin{equation}{section}% 

\newcommand{\nofigure}[1]{#1}

\newtheorem{theorem}{Theorem}[section]
\newtheorem{proposition}{Proposition}[section]
\newtheorem{lemma}{Lemma}[section]
\newtheorem{definition}{Definition}[section]

\theoremstyle{definition}
\newtheorem{remark}{Remark}[section]

\DeclareMathOperator{\supp}{supp}

\DeclareMathOperator{\WF}{WF}
\DeclareMathOperator{\dist}{dist}

\newcommand{\eps}{\varepsilon}

\newcommand{\R}{{\bf R}}
\newcommand{\Id}{\mbox{Id}}
\renewcommand{\r}[1]{(\ref{#1})}
\newcommand{\PDO}{$\Psi$DO}
\newcommand{\be}[1]{\begin{equation}\label{#1}}
\newcommand{\ee}{\end{equation}}

\renewcommand{\d}{\mathrm{d}}
\renewcommand{\i}{\mathrm{i}}
\newcommand{\bo}{\partial \Omega}

\title[Thermoacoustic tomography arising in brain imaging]{Thermoacoustic tomography arising in brain imaging}

\author[P. Stefanov]{Plamen Stefanov}
\address{Department of Mathematics, Purdue University, West Lafayette, IN 47907}
\thanks{First author partly supported by a NSF  Grant DMS-0800428 and a Simons Professorship at the MSRI}

\author[G. Uhlmann]{Gunther Uhlmann}
\address{Department of Mathematics, University of Washington, Seattle, WA 98195 and Department of Mathematics, University of California, Irvine, Irvine, CA
92697}
\thanks{Second author partly supported by NSF, a Chancellor Professorship at UC Berkeley and a Senior Clay Award}

\usepackage{graphicx}
\graphicspath{%
    {converted_graphics/}% inserted by PCTeX
    {/}% inserted by PCTeX
}
\begin{document}
\begin{abstract}
We study the mathematical model of thermoacoustic and photoacoustic tomography when the sound speed has a jump across a smooth  surface. This models the change of the sound speed in the skull when trying to image the human brain. We derive an explicit inversion formula in the form of a convergent Neumann series under the assumptions that all singularities from the support of the source reach the boundary. 
\end{abstract}

\maketitle

\section{Introduction}  
In this paper, we study the mathematical model of thermoacoustic and photoacoustic tomography for sound speed that jumps across a smooth closed surface. The physical setup is the following. A short impulse of microwaves or light is sent through  a patient's body. The cells react by emitting an acoustic signal that is being detected on a smooth surface around the patient's body. Then the problem is to recover the density of the source of the acoustic waves that can be used to recover  the absorption rate of the tissue at any point, thus creating an image, see e.g.,  \cite{Haltmeier04, Haltmeier05, Kruger03, Kruger99,XuWang06, SU-thermo}. For more detail, an extensive list of references, and the recent progress in the mathematical understanding of this problem in the case of a constant or a smooth sound speed, we refer to \cite{AgrKuchKun2008,FinchRakesh08,Hristova08, HristovaKu08,KuchmentKun08, Patch04}.

Let $\Omega\subset\R^n$ be a bounded domain with smooth boundary. 
Let $\Gamma\subset\Omega$ be a smooth closed, orientable, not necessarily connected 
surface. . 
Let the sound speed $c(x)>0$ be smooth up to $\Gamma$ with a nonzero jump across it. For $x\in\Gamma$, and a fixed orientation of $\Gamma$, we introduce the notation
\be{1.1}
c_\text{int}(x) = c\big|_{\Gamma_\text{int}}, \quad c_\text{ext}(x) = c\big|_{\Gamma_\text{ext}}
\ee
for the limits from the ``interior'' and from the ``exterior'' of $\Omega\setminus\Gamma$. 
Our assumption then is that those limits are positive as well, and 
\be{1.2}
c_\text{int}(x)\not= c_\text{ext}(x), \quad\forall x\in \Gamma.
\ee
This problem was proposed by Lihong Wang at the meeting in Banff on inverse transport and tomography in May, 2010 and it arises in brain imaging \cite{XuWang2006, YangWang2008}. In that case,  the brain is represented by some domain $\Omega_0\Subset\Omega$. Let $\Omega_1$ be another domain representing the brain and the skull, so that $\Omega_0\Subset\Omega_1\Subset\Omega$, and $\bar\Omega_1\setminus \Omega_0$ is the skull, see Figure~\ref{fig:skull2}. 
The measuring devices are then typically placed on a surface encompassing the skull, modeled by $\bo$ in our case. Then  
\[
c|_{\Omega_0}<c|_{\Omega_1\setminus \Omega_0}, \quad c|_{\Omega_1\setminus \Omega_0}> c|_{\Omega\setminus \Omega_1},
\]
with the speed jumping by about a factor of two inside the skull $\bar\Omega_1\setminus \Omega_0$. Another motivation to study this problem is to model the classical case of a smooth speed in the patient's body but account for a possible jump of the speed when the acoustic waves leave the body and enter the liquid surrounding it.

The mathematical model can then be described as follows. 
Let 
$u$ solve the problem
\begin{equation}   \label{1}
\left\{
\begin{array}{rcll}
(\partial_t^2 -c^2\Delta)u &=&0 &  \mbox{in $(0,T)\times \R^n$},\\
u\big|_{\Gamma_\text{int}}&=&  u\big|_{\Gamma_\text{ext}},\\
\frac{\partial u}{\partial \nu} \big|_{\Gamma_\text{int}}
&=& \frac{\partial u}{\partial \nu}\big|_{\Gamma_\text{ext}},\\
u|_{t=0} &=& f,\\ \quad \partial_t u|_{t=0}& =&0, 
\end{array}
\right.               
\end{equation}
where $T>0$ is fixed,  $u|_{\Gamma_\text{int,ext}}$ is the limit value (the trace) of $u$ on $\Gamma$ when taking the limit from the ``exterior'' and from the ``interior'' of $\Gamma$, respectively. We similarly define the interior/exterior normal derivatives, and $\nu$ is the exterior unit (in the Euclidean metric) normal to $\Gamma$.

Assume that $f$ is supported in $\bar\Omega$, where $\Omega\subset \R^n$ is some smooth bounded domain. The measurements are modeled by the operator
\be{1b}
\Lambda_1 f : = u|_{[0,T]\times\partial\Omega}.
\ee
The problem is to reconstruct the unknown $f$. 

We denote more generally by $\Lambda[f_1,f_2] = \Lambda_1f_1+\Lambda_2f_2$ the measurements corresponding to general Cauchy data $[f_1,f_2]$ in \r{1}. 
We will work with $f_1$, $f_2$, supported in some compact $\mathcal{K}$ in $\Omega$. In applications, this corresponds to $f_1$, $f_2$, that are not necessarily zero outside $\mathcal{K}$ but are known there. By subtracting the known part, we arrive  at the formulation that we described above. We also assume that $c=1$ on $\R^n\setminus \Omega$. We formulate the main results for Cauchy data $[f_1,0]$ for simplicity of the exposition but we do most of the preparatory work for general Cauchy data. 

The propagation of singularities for the transmission problem is well understood, at least away from possible gliding rays \cite{Hansen84, Taylor76, Petkov82, Petkov82II}. When a singularity traveling along a geodesic hits the interface $\Gamma$ transversely, there is a reflected ray carrying a singularity, that reflects at $\Gamma$ according to the usual reflection laws. If the speed on the other side is smaller, there is a transmitted (refracted) ray, as well, at an angle satisfying Snell's law, see \r{1.14}. In the opposite case, such a ray exists only if the angle with  $\Gamma$ is above some critical one, see \r{1.13c}. If that angle is smaller than the critical one, there is no transmitted singularity on the other side of $\Gamma$. This is known as a full internal reflection.  This is  what happens in the case of the skull when a ray hits the skull boundary from inside at a small enough angle, see Figure~\ref{fig:skull2}. 
 Therefore, the initial ray splits into two parts, or does not split; or hits the boundary exactly with an angle equal to the  critical one. The latter case  is more delicate, and we refer to section~\ref{sec_ans} for some discussion on that. 
Next, consider the propagation of each branch, if more than one. Each branch may split into two, etc. In the skull example, a ray coming from the interior of the skull hitting the boundary goes to a region with a smaller speed; and therefore there is always a transmitted ray, together with the reflected one. 
Then a single singularity starting at time $t=0$ until time $t=T$ in general propagates along a few branches that look like a directed graph. This is true at least under the assumption than none of those branches, including possible transmitted ones, is tangent to the boundary.

\nofigure{
\begin{figure}[t] % float placement: (h)ere, page (t)op, page (b)ottom, other (p)age
  \centering
  % file name: C:/Users/Plamen/Documents/My PCTeX Files/current/skull2.pdf
  \includegraphics[bb=0 1 800 537,width=4.09in,height=2.74in,keepaspectratio]{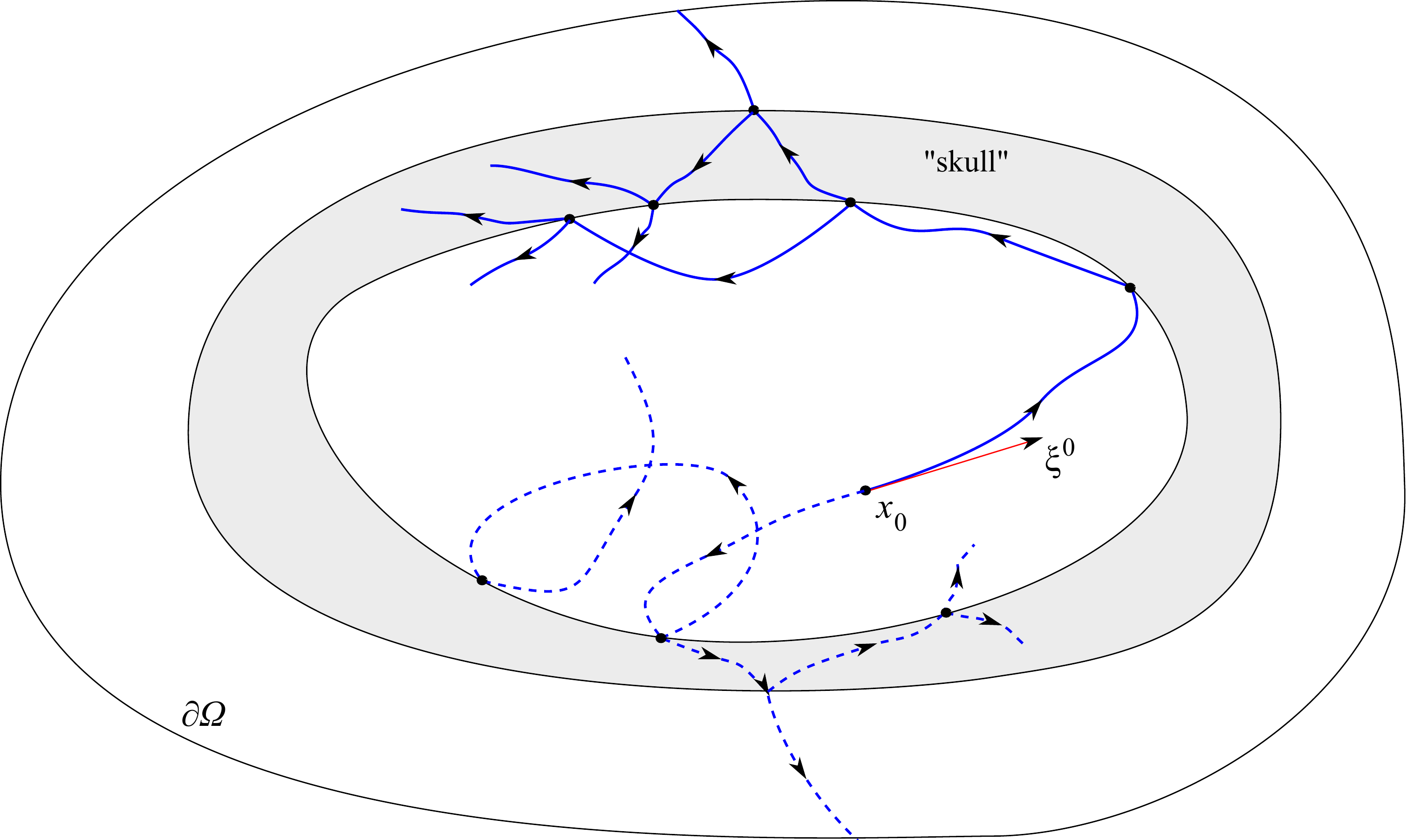}
  \caption{Propagation of singularities for the transmission problem in the ``skull'' example. The shaded region represents the ``skull'', and the speed there is higher than in the non-shaded part. The dotted curves represent the propagation of the same singularity but moving with the negative wave speed.}
  \label{fig:skull2}
\end{figure}
}

If $f_2=0$, that is the case we are interested in, singularities from $(x_0,\xi^0)$ start to propagate in the direction $\xi^0$ and in the negative one $-\xi^0$. 
If none of the branches reaches $\bo$ at time $T$ or less, a stable recovery is not possible \cite{SU-JFA}. In section~\ref{sec_whole_b}, we study the case where the initial data is supported in some compact $\mathcal{K}\subset\Omega\setminus\Gamma$ and for each $(x_0,\xi^0)\in T^*\mathcal{K}\setminus 0$, each ray through it, and through $(x_0,\xi^0)$ has a branch that reaches $\bo$ transversely at time less than  $T$. The main idea of the proof is to estimate the energy that each branch carries at high energies. If there is branching into non-tangent to the boundary rays, we show that a positive portion of the energy is transmitted, and a positive one is reflected, at high energies. As long as one of these branches reaches the boundary transversely, at a time at which measurements are still done, we can detect that singularity.  If we can do that for all singularities originating from $\mathcal{K}$, we have stability. This explains condition \r{2.1b} below.  Uniqueness follows from unique continuation results. 

Similarly to \cite{SU-thermo}, assuming  \r{2.1b}, we also get an explicit converging Neumann series formula for reconstructing $\mathbf{f}$, see  Theorem~\ref{thm2.1}. As in the case of a smooth speed considered in \cite{SU-thermo} the ``error'' operator ${K}$ in \r{2.2} is a contraction. An essential
difference in this case is that $K$ is not necessarily compact. Roughly speaking, ${K}f$ corresponds to that part of the high frequency energy that is still held in $\Omega$ until  time $T$ due to reflected  or transmitted signals that have not reached $\bo$ yet. While the first term only in \r{2.2}  will still recover all singularities of $\mathbf{f}$, it will not recover their strength, in contrast to the situation in \cite{SU-thermo}, where the speed is smooth. Thus one can expect somewhat slower convergence in this case. 

\medskip 
\noindent{\bf Acknowledgment.} We would like to thank Lihong Wang for helpful conversations on brain imaging using thermoacoustic tomography. 

\section{Main result
}\label{sec_whole_b}

Let $u$ solve the problem
\begin{equation}   \label{1'}
\left\{
\begin{array}{rcll}
(\partial_t^2 -c^2\Delta)u &=&0 &  \mbox{in $(0,T)\times \R^n$},\\
u\big|_{\Gamma_\text{int}}&=&  u\big|_{\Gamma_\text{ext}},\\
\frac{\partial u}{\partial \nu} \big|_{\Gamma_\text{int}}
&=& \frac{\partial u}{\partial \nu}\big|_{\Gamma_\text{ext}},\\
u|_{t=0} &=& f_1,\\ \quad \partial_t u|_{t=0}& =&f_2, 
\end{array}
\right.               
\end{equation}
where $T>0$ is fixed. Compared to \r{1}, we allowed $f_2$ to be non-zero but in the main result, we will take $f_2=0$ for simplicity of the exposition. 
Set $\mathbf{f}=[f_1,f_2]$. 
Assume that $f$ is supported in $\bar\Omega$, where $\Omega\subset \R^n$ is some smooth bounded domain. Set
\be{l1}
{\Lambda} \mathbf{f} : = u|_{[0,T]\times\partial\Omega}.
\ee
The trace $\Lambda \mathbf{f}$ is well defined in $C_{(0)}\big([0,T]; \; H^{1/2}(\bo)\big)$, where the subscript $(0)$ indicates that we take  the subspace of functions $h$ so that $h=0$ for $t=0$. For a discussion of other mapping properties, we refer to \cite{Isakov-book}, when $c$ has no jumps. By finite speed of propagation, one can reduce the analysis of the mapping properties of $\Lambda$ to that case. 

In the thermoacoustic model, $f_2=0$. For this reason we set
\be{l1a}
\Lambda_1 f_1 := \Lambda[f_1,0].
\ee
This notation is justified by setting  $\Lambda_{1,2}$ to be the components of $\Lambda$ (that sends vector functions to scalar functions, i.e., $\Lambda\mathbf{f} = \Lambda_1f_1+\Lambda_2f_2$. In this paper, we use boldface to denote vector functions or operators that map scalar or vector functions to vector functions.  

The standard back-projection that would serve as some kind of approximation of the actual solution is the following. We cut off smoothly $\Lambda_1f_1$ near $t=T$ to satisfy the compatibility conditions in the next step; and then we solve a backward mixed problem with boundary data the so cut $\Lambda_1f_1$; and Cauchy data $[0,0]$ at $t=T$. As in the case of a smooth speed, see \cite{Hristova08, SU-thermo}, one can show that such a back-projection would converge to $f$, as $T\to\infty$ at a rate that depends of $f$; and at least at a slow logarithmic one, if one knows a priori that $f\in H^2$, see \cite{Bellassoued03}. If $\Gamma=\bo_0$, where $\Omega_0\subset\Omega$ is strictly convex, then in the case that the speed outside $\Omega_0$ is faster than the speed inside (then there is full internal reflection), the convergence would be no faster than logarithmic,  as suggested by the result in \cite{PopovVodev_trapping}. In the opposite case, it is exponential if $n$ is odd, and polynomial when $n$ is even \cite{CardosoPV99}. Our goal in this work is to fix $T$ however.

In \cite{SU-thermo}, we proposed the following modified back-projection. 
Given $h$ that we will be chosen to be $\Lambda_1f_1$ later, let $v$ solve
\begin{equation}   \label{l2}
\left\{
\begin{array}{rcll}
(\partial_t^2 -c^2\Delta)v &=&0 &  \mbox{in $(0,T)\times \Omega$},\\
v|_{[0,T]\times\partial\Omega}&= &h,\\
v|_{t=T} &=& \varphi,\\ \quad   \partial_t v|_{t=T}& =&0, \\
\end{array}
\right.               
\end{equation}
where $\phi$ solves the elliptic boundary value problem
\be{l3}
\Delta\phi=0, \quad 
\phi|_{\partial\Omega} = h(T,\cdot).
\ee 
Then we define the following pseudo-inverse
\be{l4}
\mathbf{A} h := [v(0,\cdot),v_t(0,\cdot)] =: [A_1h,A_2h] \quad \mbox{in $\bar\Omega$}.
\ee
By  \cite{LasieckaLT}, and using finite speed of propagation \cite{Williams-transmission}, one can show that 
\[
\mathbf{A} : H^1_{(0)}([0,T]\times \bo) \to \mathcal{H}  \cong H_0^1(\Omega)\times L^2(\Omega)
\]
is a continuous map. Note that the mapping properties above allow us to apply $\mathbf{A}$ to $\Lambda\mathbf{f}$ only when $\mathbf{f}$ is compactly supported in $\Omega$ but the theorem above shows that $\mathbf{A}\Lambda$ extends continuously to the whole $\mathcal{H}(\Omega)$. The function $A_1h$ with $h=\Lambda_1f_1$ can be thought of as the first ``approximation'' of $f_1$. On the other hand, the proof Theorem~\ref{thm2.1} below shows that it is not a good approximation, see Remark~\ref{remark1}. 

To explain the idea behind this approach, let us assume for a moment that we knew the Cauchy data $[u,u_t]$ on $\{T\}\times\Omega$. Then one could simply solve the mixed problem in $[0,T]\times\Omega$ with that Cauchy data and boundary data $\Lambda_1f_1$. Then that solution at $t=0$ recovers $f$. We do not know the Cauchy data $[u,u_t]$ on $\{T\}\times\Omega$, of course, but we know the trace of $u$ (a priori in $H^1$ for $t$ fixed) on $\{T\}\times\bo$. The trace of $u_t$ does not make sense because the latter is only in $L^2$ for $t=T$. The choice of the Cauchy data in \r{l2} can then be explained by the following. Among all possible Cauchy data that belong to the ``shifted linear space'' (the  linear space $\mathcal{H}(\Omega)$ translated by a single element of the set below)
\be{C}
\left\{\mathbf{g}=[g_1,g_2]\in H^1(\Omega)\oplus L^2(\Omega); \; g_1|_{\bo}=h(T,\cdot)\right\}
\ee
we chose the one that minimizes the energy. The ``error'' then will be minimized. We refer to Lemma~\ref{lemma_D} and the proof of Theorem~\ref{thm2.1} for more details.

The discussion in the Introduction and in section~\ref{sec_ans} indicates that the singularities that we are certain to detect at $\bo$ lie in the following ``non-trapped'' set
\be{2.1a'}
\begin{split}
\mathcal{U} = &\big\{(x,\xi)\in S^*(\Omega\setminus\Gamma) ;\;  \text{there is a path  of the  ``geodesic'' issued from either}\\
& \text{ $(x,\xi)$ or $(x,-\xi)$  at $t=0$ never tangent to $\Gamma$, that is outside $\bar\Omega$ at time $t=T$} \big\}.
\end{split}
\ee
Actually, $\mathcal{U}$ is the maximal open set with the property that a singularity in $\mathcal{U}$ ``is visible'' at $[0,T]\times\bo$; and what happens at the boundary of that set, that includes for example rays tangent to $\Gamma$, will not be important for our analysis. We emphasize here that ``visible'' means that some positive fraction of the energy and high frequencies  can be detected as a singularity of the data; and of course there is a fraction that is reflected; then some trace of it may appear later on $\bo$, etc. 

One special case is the following. Take a compact set $\mathcal{K}\subset\Omega\setminus \Gamma$, and assume that 
\be{2.1b}
S^*\mathcal{K}\subset \mathcal{U}.
\ee
In other words, we require that for any $x\in \mathcal{K}$ and any unit $\xi\in S_x^*\mathcal{K}$, at least one of the multi-branched ``geodesics'' starting from $(x,\xi)$, 
and  from $(x,-\xi)$, at $t=0$ has a path that hits $\bo$ for time $t<T$ and satisfies the non-tangency assumption of \r{2.1a'}. Such a set may not even exist for some speeds $c$. 

\medskip 
\textbf{Example 1.} 
Let $\Omega_0\subset\Omega$ be two concentrating balls, and let $c$ be piecewise constant; more precisely, assume 
\[
\Omega=B(0,R), \quad \Omega_0=B(0,R_0), \quad 0<R_0<R,
\]
and let
\[
c=c_0<1 \quad \text{in $\Omega_0$}; \qquad c=1\quad \text{in $\R^n\setminus\Omega_0$}. 
\]
Then such a set $\mathcal{K}$ always exist and can be taken to be a ball with the same center and small enough radius. Indeed, the requirement then is that all rays starting from $\mathcal{K}$ hit $\Gamma$ at an angle greater than a critical one $\pi/2-\alpha_0$, see \r{1.13c}. This can be achieved by choosing $\mathcal{K}=B(0,\rho)$ with $\rho\ll R_0$. An elementary calculation shows that we need to satisfy the inequality  $\rho/R_0<\sin\alpha_0=c_0$, i.e., it is enough to choose $\rho<c_0R_0$. Then there exists $T_0$ that is easy to compute so that for $T>T_0$, \r{2.1b} holds. If $c=c_0>1$ in $\Omega_0$, then any compact $\mathcal{K}$ in $\Omega$ satisfies \r{2.1b}. In that case, there is always a transmitted ray leaving $\Omega_0$. 
\medskip  

\textbf{Example 2.} This is a simplified version of the skull model. Let $\Omega_0\subset\Omega_1 \subset \Omega$ be balls so that 
\[
\Omega=B(0,R), \quad \Omega_0=B(0,R_0),\quad \Omega_1=B(0,R_1) \quad 0<R_0<R_1<R,
\]
Assume that 
\[
c|_{\Omega_0}=c_0, \quad c|_{\Omega_1\setminus \Omega_0}=c_1, 
%\quad   c|_{\Omega\setminus \Omega_1}=c_2, 
\quad c|_{\R^n\setminus \Omega_1}=1
\]
with some constants $c_0$, $c_1$ so that $c_0<c_1$, $c_1>1$. Here, $c_0$ models the acoustic speed in the brain,  $c_1$ is the speed in the skull, and $1$ is the acoustic speed in the liquid outside the head. If for a moment we consider $\Omega_0$ and $\Omega_1$ only, we have the configuration of the previous example. If $\mathcal{K}=B(0,\rho)$ with  $\rho<(c_0/c_1)R_0$, then $\mathcal{K}$ satisfies \r{2.1b}. Now, since $c_1>1$, rays that hit $\bo_1$ always have a transmitted part outside $\Omega_1$, and therefore \r{2.1b} is still satisfied in $\Omega$.

\medskip

Let $\Pi_{\mathcal{K}}: H_D(\Omega)\to H_D(\mathcal{K})$ be the orthogonal projection of elements of the former space to the latter (considered as a subspace of $H_D(\Omega)$). It is easy to check that $\Pi_{\mathcal{K}}f= f|_{\mathcal{K}} - P_{\partial\mathcal{K}}(f|_{\partial \mathcal{K}})$, where $P_{\partial \mathcal{K}}$ is the Poisson operator of harmonic extension in $\mathcal{K}$, see also Lemma~\ref{lemma_D}. Let $\mathbf{1}_{\mathcal{K}}$ be the restriction of functions defined on $\Omega$, to $\mathcal{K}$. Then $\mathbf{\Pi}_{\mathcal{K}}\mathbf{f}= \Pi_{\mathcal{K}}f_1+\mathbf{1}_{\mathcal{K}}f_2$ is the orthogonal projection of $\mathcal{H}(\Omega)$ to $\mathcal{H}(\mathcal{K})$, similarly to \r{l3}. 

Our main result is the following.

\begin{theorem}  \label{thm2.1} 
Let $\mathcal{K}$ satisfy \r{2.1b}. Then  ${\Pi}_{\mathcal{K}}{A}_1\Lambda_1=\Id-{K}$ in $H_D(\mathcal{K)}$, with  
 $\|{K}\|_{H_D(\mathcal{K})}<1$. In particular, $\Id-{K}$ is invertible on $H_D(\mathcal{K})$, and $\Lambda_1$ restricted to $H_D(\mathcal{K})$  has an explicit left inverse of the form
\be{2.2}
{f} = \sum_{m=0}^\infty  {K}^m  {\Pi}_{\mathcal{K}} A_1 h, \quad h= \Lambda_1 f.
\ee

\end{theorem}

\begin{remark}\label{remark1}
As discussed in the Introduction, $K$  is not a compact operator as in the case of smooth sound speed. It follows from the proof of the theorem that the least upper bound of its essential spectrum (always less that $1$) corresponds to the maximal portion of the high-frequency energy that is still held  in $\Omega$ at time $t=T$.  
\end{remark}

\begin{remark}
Compared to \cite[Theorem~1]{SU-thermo}, there is a slight improvement in the theorem above; we reduce the lower bound of $T$ by a factor of $2$. Indeed, the proof also works in the case when $c$ is smooth everywhere. In \cite[Theorem~1]{SU-thermo}, we assumed $T>T(\Omega)$, where the latter is the length of the longest geodesic in $\bar\Omega$. In particular, $T(\Omega)$ is the diameter of $\bar\Omega$, if the metric $c^{-2}\d x^2$ is simple, see \cite{SU-JAMS}.  In the theorem above, $T$ has to be larger than the ``radius'' of $\bar\Omega$, see \r{2.1a'}. Therefore, if $T>T(\Omega)/2$, for any $(x,\xi)\in S\Omega$, at least one of the geodesics $\gamma_{x,\xi}$ and $\gamma_{x,-\xi}$ would exit $\bar\Omega$ at time $T$, hence \r{2.1a'} holds. Then $K$ is still a contraction but not a compact operator anymore in general, even when $c$ is smooth.  Thus one can expect a slower convergence of \r{2.2} if $T(\Omega)/2<T <T(\Omega)$.
\end{remark}

\begin{remark}
Consider the case now where $\mathcal{K}$ does not satisfy \r{2.1b}. If there is an open set of singularities that does not reach $\bo$, a stable recovery is impossible \cite{SU-JFA}. In either case however, a truncated version of the series \r{2.2} would provide an approximate parametrix that would recover the visible singularities, i.e., those in $\mathcal{U}$. By an approximate parametrix we mean a pseudo-differential operator elliptic in $\mathcal{U}$ with a principal symbol converging to $1$ in any compact in that set as the number of the terms in \r{2.2} increases. This shows that roughly speaking, if a recovery of the singularities is the primary goal, then only those in $\mathcal{U}$ can be recovered in a ``stable way'', and \r{2.2} works in that case as well, without the assumption \r{2.1b}. 
Details will be given in a forthcoming paper. 
\end{remark}

\section{Preliminaries} \label{sec_prel} 
Notice first that $c^2\Delta$ is formally self-adjoint w.r.t.\ the measure $c^{-2}\d x$.   
Given a domain $U$, and a function $\mathbf{f}=[f_1,f_2]$, define the energy
\[
E_U(\mathbf{f}) 
=  \int_U\left( |\nabla_x f_1|^2    +c^{-2}|f_2|^2 \right)\d x.
\]
Given a scalar function $u(t,x)$, we also set $\mathbf{u}=[u,u_t]$; and then 
\[
E_U(\mathbf{u}(t)) = \int_U\left( |\nabla_x u|^2   +c^{-2}|u_t|^2 \right)\d x.
\] 
Here and below we use the notation $\mathbf{u}(t)=\mathbf{u}(t,\cdot)$. 
We define the energy space $\mathcal{H}(U)$ as the  completion of $C_0^\infty(U)\times  C_0^\infty(U)$ under the energy norm
\[
\|\mathbf{f}\|^2_{\mathcal{H}(U)} = E_U(\mathbf{f}) = \int_U\left( |\nabla_x f_1|^2    +c^{-2}|f_2|^2 \right)\d x.
\]
In particular, we define the space $H_{D}(U)$ to be the completion of $C_0^\infty(U)$ under the Dirichlet norm
\be{2.0H}
\|f\|_{H_{D}(U)}^2= \int_U |\nabla_x u|^2 \,\d x.
\ee
The energy norm is topologically equivalent to the same norm with $c=1$. 
It is easy to see that $H_{D}(U)\subset H^1(U)$, if $U$ is bounded with smooth boundary, therefore, $H_{D}(U)$ is topologically equivalent to $H_0^1(U)$. If $U=\R^n$, this is true for $n\ge3$ only \cite{LP}.  By finite speed of propagation, the solution with compactly supported Cauchy data always stays in $H^1$ even when $n=2$.  
This norm defines the energy space 
\[
\mathcal{H}(U) = H_D(U)\oplus L^2(U).
\] 
Here and below, $L^2(U) = L^2(U; \; c^{-2}\d x)$. Note also that 
\be{Pf}
\|f\|^2_{H_D(U)} = (-c^2\Delta f,f)_{L^2}.
\ee
The wave equation then can be written down as the system
\be{s1}
\mathbf{u}_t= \mathbf{P}\mathbf{u}, \quad \mathbf{P} = \begin{pmatrix} 0&I\\c^2\Delta&0 \end{pmatrix},%\quad P := c^2\Delta,
\ee
where $\mathbf{u}=[u,u_t]$ belongs to the energy space $\mathcal{H}$. 

\begin{proposition}
Let $U$ be an open subset of $\R^n$ with smooth boundary. Then the operator $\mathbf{P}$  with domain   $D(\mathbf{P})=\{\mathbf{f}\in \mathcal{H}(U); \; f_2\in H_D(U), \; \Delta f_1\in L^2(U) \}$ is a unbounded skew-selfadjoint operator on $\mathcal{H}(U)$ if $c\in L^\infty$, and $c^{-1}\in L^\infty$. 
\end{proposition}

\begin{proof}
Clearly, $\mathbf{P}$ is correctly defined on $D(\mathbf{P})$, and the latter is dense in $\mathcal{H}(U)$.  We will show first that $\mathbf{P}^*\supset -\mathbf{P}$, i.e., that $D(\mathbf{P}^*)\supset D(\mathbf{P})$ and $\mathbf{P}^*=-\mathbf{P}$ on $D(\mathbf{P})$. Let $\mathbf{f}$, $\mathbf{g}$ be in $D(\mathbf{P})$. Then 
\be{0.1}
(\mathbf{Pf},\mathbf{g})_{\mathcal{H}(U)} = \int_U \left(  \nabla f_2 \cdot\nabla \bar g_1 -   \nabla f_1 \cdot\nabla \bar g_2 \right)\, \d x.
\ee
This follows by integration by parts that can be easily justified by approximating $\mathbf{g}$ with smooth functions, see also \cite[Theorem~V.1.2]{LP}. Replace $\mathbf{f}$ by 
and $\mathbf{\bar f}$,  and $\mathbf{g}$ by $\mathbf{\bar g}$ to get
\[
(\mathbf{Pf},\mathbf{g})_{\mathcal{H}(U)} = -(\mathbf{f},\mathbf{Pg})_{\mathcal{H}(U)} , \quad \forall \mathbf{f}, \mathbf{g}\in D(\mathbf{P}).
\]
We will show next that $\mathbf{P}^*\subset -\mathbf{P}$. For that, in view of what we already proved, it is enough to show that $D(\mathbf{P}^*)=D(\mathbf{P})$. By definition, $\mathbf{g}\in D(\mathbf{P}^*)$ if and only if there exists $\mathbf{h}\in \mathcal{H}(U)$ so that
\[
(\mathbf{Pf},\mathbf{g})_{\mathcal{H}(U)} = (\mathbf{f},\mathbf{h})_{\mathcal{H}(U)}
,\quad \forall \mathbf{f}\in D(\mathbf{P}),
\]
see e.g., \cite[\S 8.1]{Reed-Simon1}. This equality can also be written as
\[
\int_U \left(  \nabla f_2 \cdot\nabla \bar g_1 +  (\Delta f_1) \bar g_2 \right)\, \d x = \int_U \left(\nabla f_1\cdot\nabla\bar h_1  +c^{-2} f_2\bar h_2            \right)\,\d x.
\]
Note that the l.h.s.\ is independent of $c$. Set $\mathbf{h}^\sharp=[h_1,c^{-2}h_2]$. It belongs to $\mathcal{H}$ if and only if $\mathcal{h}$ does. 
We therefore reduced the problem $D(\mathbf{P}^*)=D(\mathbf{P})$ to its partial case when $c=1$. In that case, this statement is well known, see e.g., \cite[Theorem~V.1.2]{LP}. In other words, such $\mathbf{h}^\sharp\in \mathcal{H}_1(U)$   exists, and equals $-\mathbf{P}_1\mathbf{g}$, where $\mathcal{H}_1$, $\mathbf{P}_1$ correspond to $c=1$, if and only if  $\mathbf{g}\in D(\mathbf{P}_1)=D(\mathbf{P})$. Then we also get that $\mathbf{h}\in \mathcal{H}(U)$ exists (and, of course, equals $-\mathbf{Pf}$). 
\end{proof}

\begin{remark}
If $U$ is bounded then $D(\mathbf{P})$ is topologically equivalent to $H^2(U)\oplus H^1(U)$. If $U$ is unbounded but $\supp \mathbf{f}$ is a compact in $U$, then any $\mathbf{f}\in D(\mathbf{P})$ belongs to that space as well. The case of non-compactly supported $\mathbf{f}$ is more delicate, see \cite{LP}. In this paper, we will deal with either $U=\R^n$ or $U=\Omega$; and in the former case, $\mathbf{f}$ is compactly supported. In the latter case, the definition of $H_D(U)$ reflects Dirichlet boundary conditions on $\bo$.
\end{remark}

By Stone's theorem, $\mathbf{P}$ is a generator of a strongly continuous group $e^{t\mathbf{P}}$ of unitary operators on $\mathcal{H}(U)$ that leaves $D(\mathbf{P})$ invariant. For any $\mathbf{f}\in D( \mathbf{P})$, we have $\frac{\d}{\d t}e^{t\mathbf{P}}\mathbf{f}= \mathbf{P}e^{t\mathbf{P}}\mathbf{f}$. In particular, for the first component $u(t,x)$ of $e^{t\mathbf{P}}\mathbf{f}$, we have that $u_{tt}(t,\cdot)$ is a continuous function of $t$ with values in $L^2(U)$; the same holds for $c^2(\cdot)\Delta u(t,\cdot)$, and $u_{tt}=c^2\Delta u$. 

\begin{remark}
An alternative proof of the proposition is to use well-posedness results in \cite{Williams-transmission}, for example, to show that the solution is given by a strongly continuous unitary group; and deduce from there that its generator  must be skew-selfadjoint.
\end{remark}

For $\mathbf{f}\in D(\mathbf{P})$, the following transmission conditions hold for $f_1$
\be{1.3}
f_1\big|_{\Gamma_\text{int}}=  f_1\big|_{\Gamma_\text{ext}},\quad 
\frac{\partial f_1}{\partial \nu} \Big|_{\Gamma_\text{int}}
= \frac{\partial f_1}{\partial \nu}\Big|_{\Gamma_\text{ext}},
\ee
that just reflects the fact that $f_1$ is in $H^2$ locally. 

\section{Geometric Optics}\label{sec_ans}
\subsection{A parametrix of the Cauchy problem with data at $t=0$} 
We start with a standard geometric optics construction. We assume first that $c$ is smooth in the region where we construct the parametrix. 
 Fix $(x_0,\xi^0)\in T^*\Omega\setminus 0$. We identify vectors and covectors away from $\Gamma$ by the metric $c^{-2}\d x^2$. 
Let $u$ be a solution of the wave equation with initial data $[u,u_t]=[f_1,f_2]$, with $\mathbf{f}= [f_1,f_2]$ having a wave front set in some small conic neighborhood of $  (x_0,\xi^0)$. 

In a neighborhood of $(0,x_0)$, the solution to \r{1'} is given by
\be{o1}
u(t,x) =  (2\pi)^{-n} \sum_{\sigma=\pm}\int e^{\i\phi_\sigma(t,x,\xi)}\left( a_{1,\sigma}(x,\xi,t) \hat f_1(\xi)+  a_{2,\sigma}(x,\xi,t) \hat f_2(\xi)\right)\, \d \xi,
\ee
modulo smooth terms. Here the phase functions $\phi_\pm$ are positively homogeneous of order $1$ in $\xi$ and solve the eikonal equations
\be{o2}
\pm\partial_t\phi_\pm + c(x)|\nabla_x\phi_\pm|=0,\quad \phi_\pm|_{t=0}=x\cdot\xi,
\ee
while $a_{j,\pm}$ are classical amplitudes of order $0$ 
solving  the corresponding transport equations along bicharacteristics issued from some conic neighborhood of  $(x,\pm \xi)$ for all $(x,\xi)\in \WF(\mathbf{f})$, see \cite[p.~128]{Duistermaat} or \cite[eqn.~(VI.1.50)]{Treves2}. In particular, $a_{j,\pm}$ satisfy
\[
a_{1,+} +a_{1,-}=1,\quad a_{2,+} +a_{2,-}=0        \quad \mbox{for $t=0$}.
\]
Since $\partial_t\phi_\pm=\mp c(x)|\xi|$ for $t=0$, and $u_t=f_2$ for $t=0$, we also see that for the principal terms $a_{j,\pm}^{(0)}$ of $a_{j,\pm} \sim\sum_{k\ge0} a_{j,\pm}^{(-k)}$  we have
\[
a_{1,+}^{(0)}=a_{1,-}^{(0)}, \quad c(x)|\xi|(a_{2,-}^{(0)}-a_{2,+}^{(0)})=1\quad
\mbox{for $t=0$}.
\]
Therefore, 
\be{o2a}
a_{1,+}^{(0)}=a_{1,-}^{(0)}=1/2,\quad  a_{2,+}^{(0)}=-a_{2,-}^{(0)}= \frac{1}{2c(x)|\xi|}\quad\text{for $t=0$}. 
\ee
Note that if $c=1$, then $\phi_\pm = x\cdot\xi \mp t|\xi| $, and $a_{1,+}\equiv a_{1,-}=1/2$, $a_{2,+}\equiv -a_{2,-}=(2|\xi|)^{-1}$.  The principal terms $a_{j,\pm}^{(0)}$ satisfy the homogeneous transport equations
\be{tr}
\left( (\partial_t \phi_\pm) \partial_t - c^{2} \nabla_x\phi_\pm\cdot \nabla_x+C_\pm \right)a_{j,\pm}^{(0)}=0
\ee
with initial conditions given by \r{o2a}, 
where %$C_\pm$ depend on the coefficients of $P$ and on $\phi_\pm$, see \cite[eqn.~(VI.1.49)]{Treves2}. In the particular case we study,
\[
2C_\pm = (\partial_t^2-c^2\Delta)\phi_\pm,
\]
see also \cite[eqn.~(VI.1.49)]{Treves2}.

By the stationary phase method, for $\sigma=+$, singularities starting from $(x,\xi)\in \WF(\mathbf{f})$ propagate along geodesics in the phase space  issued from $(x,\xi)$, i.e., they stay on the curve $(\gamma_{x,\xi}(t),\dot\gamma_{x,\xi}(t) )$; and from $(x,-\xi)$, for $\sigma=-$, i.e., they stay on the curve $(\gamma_{x,-\xi}(t), \dot\gamma_{x,-\xi}(t))$. This is consistent with the general propagation of singularities theory for the wave equation because the principal symbol  $-\tau^2+c^{2}|\xi|^2$ of the wave operator has two roots $\tau = \pm c|\xi|$.

The construction is valid as long as the eikonal equations are solvable, i.e., for $|t|$ small enough, at least. 
We call the corresponding parametrices in \r{o1}  $u_+$ and $u_-$.  It is well known that $u_++u_-$ represents the solution of \r{1'} up to a smoothing operator applied to $\mathbf{f}$.

\subsection{Projections to the positive and the negative wave speeds} The zeros of the principal symbol of the wave operator, in regions where $c$ is smooth,  are given by $\tau= \pm c(x)|\xi|$, that we call wave speeds. We constructed above parametrices $u_\pm$ for the corresponding solutions. We will present here a functional analysis point of view that allows us to project the initial data $\mathbf{f}$ to data $\mathbf{\Pi}_\pm\mathbf{f}$, so that, up to smoothing operators, $u_{\pm}$ corresponds to initial data $\mathbf{\Pi}_\pm\mathbf{f}$.

Assume that $c(x)$ is extended form the maximal connected component of $\R^n\setminus\Gamma$ containing $x_0$ to the whole $\R^n$ in a smooth way so that $0<1/C\le c(x)\le C$. 
%Since we are only going to need this constriction for $\mathnf{U}(t) = e^{t\mathbf{P}\mathbf{f}$ for $\mathbf{f}$ supported in $\Omega\setminus \Gamma$ and $|t|\ll1$ so that the support of, the way we choose this extension is not important. 
 Let 
\be{42.1}
Q=(-c^2\Delta)^{1/2},
\ee
where the operator in the parentheses is the natural  self-adjoint extension of $-c^2\Delta$ to $L^2(\R^n, c^{-2}\d x)$, and the square root exists by the functional calculus. Moreover, $Q$ is an elliptic \PDO\ of order $1$ in any open set; and let $Q^{-1}$ denote a fixed parametrix. 

It is well known that the solution to \r{1'} can be written as
\be{42.2}
u = \cos(tQ)f_1 + \frac{\sin(tQ)}{Q}f_2,
\ee
and the latter operator is defined by the functional calculus as $\phi(t,Q)$ with $\phi(t,\lambda)=\sin(t\cdot)/\cdot\in C^\infty$. Based on that, we can write
\be{42.3}
e^{t\mathbf{P}}= e^{\i tQ}\mathbf{\Pi}_+ + e^{-\i tQ}\mathbf{\Pi}_-,
\ee 
where 
\be{42.4}
\mathbf{\Pi}_+ = 
\frac12 \begin{pmatrix}  1   &-\i Q^{-1} \\\i Q&1 \end{pmatrix},\quad 
\mathbf{\Pi}_- = \frac12 \begin{pmatrix}  1   &\i Q^{-1} \\ -\i Q&1 \end{pmatrix}.
\ee

It is straightforward to see that $\mathbf{\Pi}\pm$ are orthogonal projections in $\mathcal{H}$, up to errors of smoothing type. Then given $\mathbf{f}\in \mathcal{H}$ supported on $\Omega$, one has $\mathbf{u}_\pm = e^{t\mathbf{P}}\mathbf{f}_\pm$, with $\mathbf{f}_\pm := \mathbf{\Pi}_\pm \mathbf{f}$.

\subsection{Analysis at the boundary} We will analyze  what happens when the geodesic   $(x_0,\xi^0)$ issued from $(x_0,\xi^0)$, $x_0\not\in \Gamma$,  hits $\Gamma$ for  first time, under some assumptions.  Let the open sets $\Omega_\text{int}$, $\Omega_\text{ext}$,   be the ``interior''  and the ``exterior'' part of $\Omega$ near $x_0$, according to the orientation of $\Gamma$. They only need to be defined near $x_0$. Let us assume that this geodesic hits $\Gamma$ from $\Omega_\text{int}$. 
 We will construct here a microlocal representation of the reflected and the transmitted waves near the boundary. 
 
 Extend $c|_{\Omega_\text{int}}$ in a smooth way in a small neighborhood on the other side of $\Gamma$, and let  $u_+$ be the solution described above, defined in some neighborhood of that geodesic segment. Since we are only going to use $u_+$ in the microlocal construction described below, and we will need only the trace of $u_+$ on $\R_+\times\Gamma$ near the first contact of the bicharacteristic from $(x_0,\xi^0)$ with $\Gamma$, the particular extension of $c$ would not affect the microlocal expansion but may affect the smoothing part. 

Set
\be{1.23h}
h := u_+|_{\R\times\Gamma}.
\ee
Let $(t_1,x_1)\in \R_+\times\Gamma$ be the point where the geodesic from $\gamma_{x_0,\xi^0}$ hits $\Gamma$ for the first time, see Figure~\ref{fig:skull1}. We assume that such $t_1$ exists. 
Let $\xi^1$ be the tangent covector to that geodesic at $(t_1,x_1)$. Assume that $\xi^0$ is unit covector in the metric $c^{-2}\d x$, then so is $\xi^1$ (in the metric  $c_\text{int}^{-2}\d x$), i.e., $c_\text{int}|\xi|=1$, where $|\xi|$ is the Euclidean norm. Assume that $\xi^1$ is transversal to $\Gamma$. In view of condition \r{2.1b}, this is the case that we need to study.

Standard microlocal arguments show, see \cite[Proposition~3]{SU-thermo} for details, that the map $[f_1,f_2]\mapsto h$ is an elliptic  Fourier Integral Operator (FIO)  with a canonical relation that is locally a canonical graph described in \cite[Proposition~3]{SU-thermo}. That diffeomorphism maps $(x_0,\xi^0)$ into $(t_1,x_1,1, (\xi^1)')$, where the prime stands for the tangential projection onto $T^*\Gamma$; and that maps extends as a positively homogeneous one of order one w.r.t.\ the dual variable. In particular, the dual variable $\tau$ to $t$ stays positive. 
In fact,  $\WF (u)$ is in the characteristic set $\tau^2-c^2(x)|\xi|^2=0$,  and $(x,\xi)$ belongs to some small neighborhood of $(x_1,\xi^1)$. 
The wave front set $\WF(h)$ is given by $(x,\xi')\in T^*\Gamma$, $(x,\xi)\in \WF(u)$, where $\xi'$ is the tangential projection of $\xi$ to the boundary. Then $(t,x,\tau,\xi')$ is the image of some $(\tilde x,\tilde\xi)$  close to $(x_0,\xi^0)$ under the canonical map above. Here $(\tilde x,\tilde\xi)$ is such that the $x$-projection $x(s)$ of the bicharacteristic  from it hits $\Gamma$ for the first time at time for the value of $s$ given by $sc(\tilde x)=t$.  Since $\tau^2-c_\text{int}^2(x)|\xi|^2=0$, for the projection $\xi'$ we have $\tau^2-c_\text{int}^2(x)|\xi'|^2>0$, where $(x,\xi')\in T^*\Gamma$, and $|\xi'|$ is the norm of the covector $\xi'$ in the metric on $\Gamma$ induced by the Euclidean one. 

The microlocal regions of $T^*(\R\times \Gamma)\ni (t,x,\tau,\xi')$ with respect to the sound speed $c_\text{int}$, i.e., in $\bar \Omega_\text{int}$,  
 are defined as follows: 
\begin{itemize}
  \item[]  \textit{hyperbolic region}: $c_\text{int}(x)|\xi'|<\tau$,
  \item[] \textit{glancing manifold}:  $c_\text{int}(x)|\xi'|=\tau$,
  \item[]  \textit{elliptic region}:  $c_\text{int}(x)|\xi'|>\tau$.
\end{itemize}
One has a similar classification of $T^*\Gamma$ with respect to the sound speed $c_\text{ext}$. A ray that hits $\Gamma$ transversely, coming from $\Omega_\text{int}$, has a tangential projection on $T^*(\R\times \Gamma)$ in the hyperbolic region relative to $c_\text{int}$. If  $c_\text{int}<c_\text{ext}$, that projection may belong to any of the three microlocal regions w.r.t.\ the speed $c_\text{int}$. If $c_\text{int}>c_\text{ext}$, then that projection is always in the hyperbolic region for $c_\text{ext}$. When we have a ray that hits $\Gamma$ from $\Omega_\text{ext}$, then those two cases are reversed. 

\subsection{The reflected and the transmitted waves}\label{sec_r_t}
We will analyze the case where $(\xi^1)'$ belongs to the hyperbolic region with respect to both $c_\text{int}$ and $c_\text{ext}$, i.e., we will work with $\xi'$ in a neighborhood of $(\xi^1)'$ satisfying 
\be{1.17c}
c_\text{int}^{-2}\tau^2- |\xi'|^2 >0, \quad c_\text{ext}^{-2}\tau^2- |\xi'|^2 >0.
\ee
The analysis also applies to the case of a ray coming from $\Omega_\text{ext}$, under the same assumption. 
We will confirm  below in this setting the well known fact  that under that condition, such a ray splits into a reflected  ray with the same tangential component of the velocity that returns to the interior  $\Omega_\text{int}$, and a transmitted one, again with the same tangential component of the velocity, that propagates in  $\Omega_\text{ext}$. We will also compute the amplitudes and the energy at high frequencies of the corresponding asymptotic solutions.

Choose local coordinates on $\Gamma$ that we denote by $x'$, and a normal coordinate $x^n$ to $\Gamma$ so that $x^n>0$ in $\Omega_\text{ext}$, and $|x^n|$ is the Euclidean distance to $\Gamma$; then $x=(x',x^n)$. 
We will express the solution $u_+$ in $\R\times\bar\Omega_\text{int}$ that we defined above, as well as a reflected solution $u_R$ in the same set; and a transmitted one $u_T$ in $\R\times\bar\Omega_\text{ext}$, up to smoothing terms in the form
\be{1.24}
u_\sigma = (2\pi)^{-n}\int e^{\i  
\varphi_\sigma(t,x,\tau,\xi')}b_\sigma(t,x,\tau,\xi')\hat h(\tau,\xi')\, \d\tau\,\d\xi',\quad
\sigma = +,R,T,
\ee
where $\hat h:= \int_{\R\times\R^{n-1}} e^{-\i (- t\tau+x'\cdot\xi') }h(t,x')\d t\, \d x'$. We chose to alter the sign of $\tau$ so that if $c=1$, then the phase function in \r{1.24} would equal $\varphi_+$, i.e., then $\varphi_+ = - t\tau+x\cdot\xi$.  
The three phase functions $\varphi_+$, $\varphi_R$, $\varphi_T$ solve  the eikonal equation %below with a choice of sign that we will describe below
\be{1.24a}
\partial_t\varphi_\sigma +c(x)|\nabla_x\varphi_\sigma|=0,\quad \varphi_\sigma|_{x^n=0}=-t\tau + x'\cdot\xi'.
\ee
The right choice of the sign in front of $\partial_t\varphi_+$, see \r{o2}, is the  positive   one because $\partial_t\varphi_+=-\tau<0$ for $x^n=0$, and that derivative must remain negative near the boundary as well. We see below that $\varphi_{R,T}$ have the same boundary values on $x^n=0$, therefore they satisfy the same eikonal equation, with the same choice of the sign.

Let now $h$ be a compactly supported distribution on $\R\times\Gamma$ with $\WF(h)$ in a small conic neighborhood of $(t_1,x_1,1,(\xi^1)')$. We will take $h$ as in \r{1.23h} eventually, with $u_+$ the solution corresponding to initial data $\mathbf{f}$ at $t=0$ but in what follows, $h$ is arbitrary as long as $\WF(h)$ has that property, and $u_+$ is determined through $h$. 
We now look for a parametrix 
\be{1.24p}
\tilde u=u_++u_R+u_T 
\ee 
near $(t_1,x_1)$ with $u_+$, $u_R$, $u_T$ of the type \r{1.24}, satisfying the wave equation and  \r{1.23h}.  We use the notation for $u_+$ now for a parametrix in $\Omega_\text{int}$ having singularities that come from the past and hit $\Gamma$; i.e.,  for an outgoing solution. The subscript $+$ is there to remind us that this is related to the positive sound speed $c(x)|\xi|$. Next, 
$u_R$ is a solution with singularities that are obtained form those of $u_+$ by reflection; they propagate back to $\Omega_\text{int}$. It is an outgoing solution in  $\Omega_\text{int}$. And finally,  $u_T$ is a solution in  $\Omega_\text{ext}$ with singularities that go away from $\Gamma $ as time increases; hence it is outgoing there.   To satisfy the first  transmission condition in \r{1.3}, we need to have
\be{1.4}
\varphi_T =\varphi_R = \varphi_+ = -t\tau+x\cdot\xi'\quad \text{for $x^n=0$},
\ee
that explains the same boundary condition in \r{1.24a}, and
\be{1.5}
1+b_R= b_T \quad \text{for $x^n=0$}.
\ee
In particular, for the leading terms of the amplitudes we get
\be{1.5a}
b_T^{(0)} -b_R^{(0)} = 1 \quad \text{for $x^n=0$}.
\ee
To satisfy the second transmission condition, we require
\be{1.6}
\i \frac{\partial \varphi_+}{\partial x^n}   + \frac{\partial b_+}{\partial x^n}   
+ \i \frac{\partial \varphi_R}{\partial x^n}   b_R + \frac{\partial b_R}{\partial x^n} = 
\i \frac{\partial \varphi_T}{\partial x^n}   b_T + \frac{\partial b_T}{\partial x^n} \quad \text{for $x^n=0$}.
\ee
Expanding this in a series of homogeneous in $(\tau,\xi)$ terms, we get series of initial conditions for the transport equations that follow. Comparing the leading order terms only, we get
\be{1.7}
  \frac{\partial \varphi_T}{\partial x^n}   b_T^{(0)}
-  \frac{\partial \varphi_R}{\partial x^n}   b_R^{(0)}
 = 
 \frac{\partial \varphi_+}{\partial x^n}     \quad \text{for $x^n=0$}.
\ee
The linear system \r{1.5a}, \r{1.7} for $b_R^{(0)}|_{x^n=0}$, $b_T^{(0)}|_{x^n=0}$ has determinant
\be{1.8}
-\left( \frac{\partial \varphi_T}{\partial x^n} - \frac{\partial \varphi_R}{\partial x^n}\right)\bigg|_{x^n=0}.
\ee
Provided  that this determinant is non-zero near $x_1$, we can solve for $b_R^{(0)}|_{x^n=0}$, $b_T^{(0)}|_{x^n=0}$. Moreover, the determination of each subsequent term $b_R^{(-j)}|_{x^n=0}$, $b_T^{(-j)}|_{x^n=0}$ in the asymptotic expansion of $b_R|_{x^n=0}$, $b_T|_{x^n=0}$ can be found by \r{1.6} by solving a linear system with the same (non-zero) determinant.

\subsection{Solving the eikonal equations} 
As it is well known, the eikonal equation \r{1.24a} in any  fixed side of $\R\times \Gamma$, near $(t_1,x_1)$, has two solutions. They are determined by a choice of the sign of the normal derivative on  $\R\times\Gamma$  and the  boundary condition. We will make the choice of the signs according to the desired properties for the singularities of $u_+$, $u_R$, $u_T$. 
Let $\nabla_{x'}$ denote the tangential gradient on $\Gamma$. By \r{1.4}, 
\be{1.10}
\nabla_{x'}\varphi_T  = 
\nabla_{x'}\varphi_R = \nabla_{x'}\varphi_+=\xi', \quad  
\partial_t \varphi_T = \partial_t \varphi_R = \partial_t\varphi_+=-\tau 
\quad \text{for $x^n=0$}.
\ee
 Using the eikonal equation \r{1.24a} and the boundary condition there, we get
\be{1.25}
\frac{\partial\varphi_+}{\partial t}= -\tau, \quad \frac{\partial\varphi_+}{\partial x^n}= 
\sqrt{c_\text{int}^{-2}\tau^2-|\xi'|^2  } 
\quad \text{for $x^n=0$}. 
\ee
We made a sign choice for the square root here based on the required property of $u_+$ described above. 
This shows in particular, that the map $h\mapsto \partial u_+/\partial t$ (that is just $\d/\d t$), and the interior incoming   Dirichlet to Neumann map 
\[
N_{\text{int,in}} : h\mapsto \frac{\partial u_+}{\partial \nu}\Big|_{\R\times\Gamma}
\]
are locally \PDO s of order $1$ with principal symbols given by $-\i\tau$, and 
\be{1.25N}
\sigma_p(N_{\text{int,in}}) = \i \frac{\partial\varphi_+}{\partial x^n}
= \i \sqrt{c_\text{int}^{-2}\tau^2-|\xi'|^2  }.
\ee
The notion ``interior incoming'' is related to the fact that locally, near $(t_1,x_1)$, we are solving a mixed problem in $\R\times \Omega_\text{int}$ with lateral boundary value $h$ and zero Cauchy data  for $t\gg0$.

\nofigure{
\begin{figure}[t] % float placement: (h)ere, page (t)op, page (b)ottom, other (p)age
  \centering
  % file name: C:/Users/Plamen/Documents/My PCTeX Files/current/skull1.pdf
  \includegraphics[bb=0 0 1103 594,width=6.21in,height=3.35in,keepaspectratio]{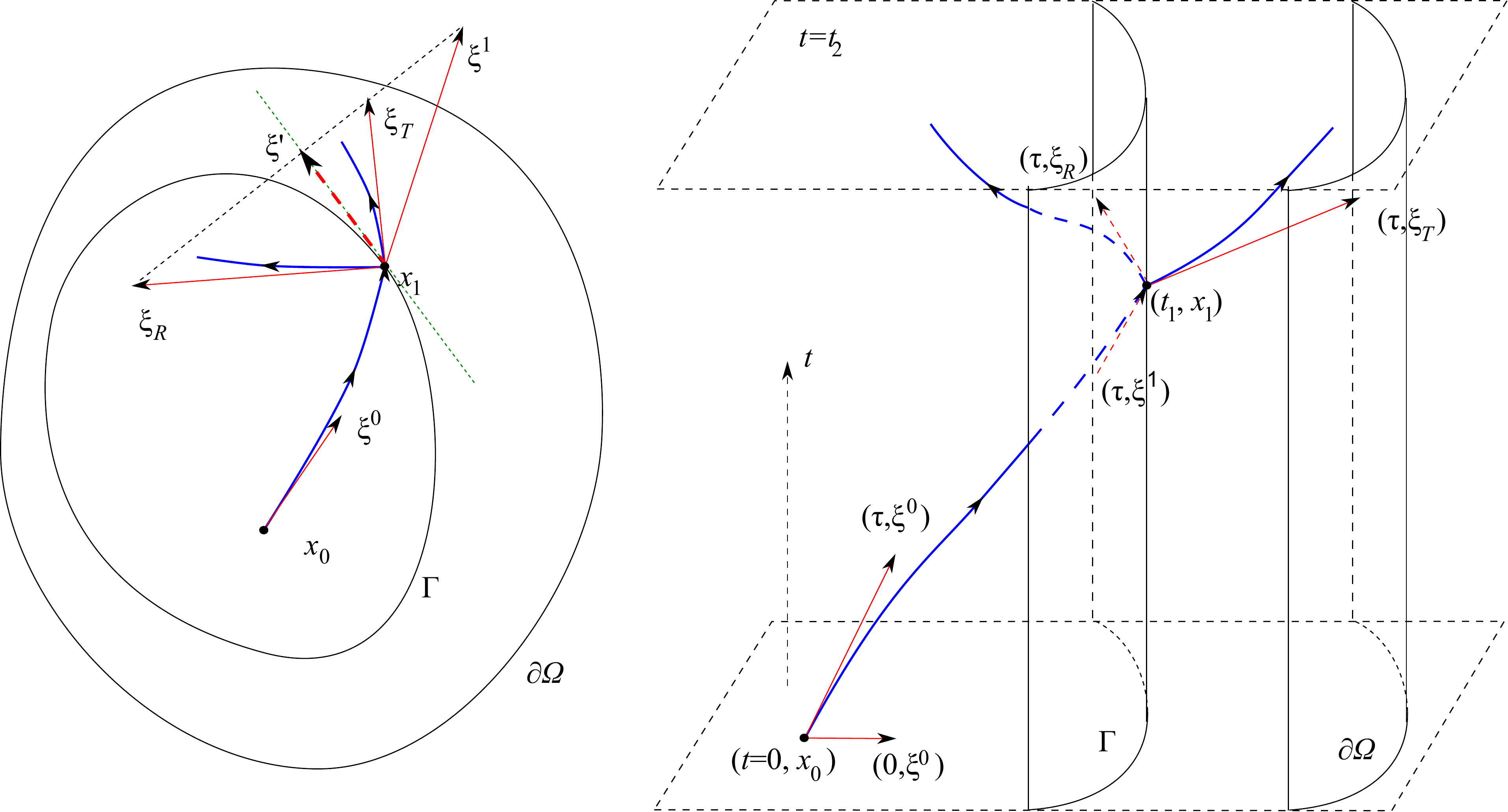}
  \caption{}
  \label{fig:skull1}
\end{figure}
}

Consider $\varphi_R$ next. The reflected phase $\varphi_R$ solves the same eikonal equation, with the same boundary condition, as $\varphi_+$. By the eikonal equation \r{1.24a}, we must have
\be{1.10a}
\frac{\partial \varphi_R}{\partial x^n}=  \pm \frac{\partial \varphi_+}{\partial x^n} \quad \text{for $x^n=0$}.
\ee
The ``$+$'' choice will give us the solution $\varphi_+$ for $\varphi_R$. We chose the negative sign, that uniquely determines a solution locally, that we call $\varphi_R$, i.e.,
\be{1.10f}
\frac{\partial \varphi_R}{\partial x^n}=  - \frac{\partial \varphi_+}{\partial x^n} \quad \text{for $x^n=0$}.
\ee
Therefore, $\nabla \varphi_R$ on the boundary is obtained from $\nabla \varphi_+$ by inverting the sign of the normal derivative. This corresponds to the usual law of reflection. 
 Therefore, 
\be{1.26}
\frac{\partial\varphi_R}{\partial t}= -\tau, \quad \frac{\partial\varphi_R}{\partial x^n}= 
-\sqrt{c_\text{int}^{-2}\tau^2- |\xi'|^2  } \quad \text{for $x^n=0$}. 
\ee
In particular,  $\partial u_R/\partial x^n|_{\R\times\Gamma}$ can be obtained from $u_R|_{\R\times\Gamma}$, that we still need to determine,  via the interior outgoing  Dirichlet-to-Neumann map
\[
N_{\text{int,out}}: u_R\Big|_{\R\times\Gamma}\quad \longmapsto \quad \frac{\partial u_R}{\partial x^n}\Big|_{\R\times\Gamma}
\]
that is locally  a first order \PDO\ with principal symbol 
\be{1.26N}
\sigma_p(N_{\text{int,out}}) =\i \frac{\partial\varphi_R}{\partial t}  =-\i \sqrt{c_\text{int}^{-2}\tau^2-|\xi'|^2  }.
\ee

To construct $\varphi_T$, we work in $\bar\Omega_\text{ext}$. We define $\varphi_T$ as the solution of \r{1.24a} with the following choice of  a normal derivative. This time $\varphi_T$ and $\varphi_+$ solve the eikonal equation at different sides of $\Gamma$, and $c$ has a jump at $\Gamma$.   By \r{1.10}, 
\be{1.11}
c_\text{ext}^2\left( |\xi'|^2   +  \Big|\frac{\partial \varphi_T}{\partial x^n}\Big|^2  \right)=\tau^2\quad \text{for $x^n=0$}.
\ee
We solve this equation for $|\partial\varphi_T/\partial x^n|^2$. Under the assumption \r{1.17c}, this solution is positive, therefore we can solve for $\partial\varphi_T/\partial x^n$ to get
\be{1.27}
 \frac{\partial\varphi_T}{\partial x^n}  = 
\quad  \sqrt{c_\text{ext}^{-2}\tau^2- |\xi'|^2 }\quad 
\text{for $x^n=0$}. 
\ee
The positive sign of the square root is determined by the requirement the singularity to be outgoing. In particular, we get that 
the exterior outgoing Dirichlet-to Neumann map
\[
N_{\text{ext,out}}: u_T\Big|_{\R\times\Gamma} \longmapsto \frac{\partial u_T}{\partial x^n}\Big|_{\R\times\Gamma}
\]
has principal symbol 
\be{1.26N2}
\sigma_p(N_{\text{ext,out}}) =\i  \frac{\partial\varphi_T}{\partial x^n}  = \i \sqrt{c_\text{ext}^{-2}\tau^2-|\xi'|^2  }.
\ee
%and again, we have an explicit formula for $b_T^{(0)}$. 

For future reference, we note that the following inequality holds
\be{1.17_0}
0\le \frac{\partial\varphi_T}{\partial x^n} \le \gamma\frac{\partial\varphi_+}{\partial x^n}; \quad
\gamma := \max_{\Gamma}\frac{c_\text{int}}{c_\text{ext}}<1.
\ee

\subsection{Amplitude and Energy Calculations} By \r{1.26}, \r{1.27}, the determinant \r{1.8} is negative. Solving \r{1.5a} and \r{1.7} then yields
\be{1.16}
b_T^{(0)} = \frac{2\partial\varphi_+/\partial x^n}{\partial\varphi_+/\partial x^n +\partial\varphi_T/\partial x^n  },\quad
b_R^{(0)} = \frac{\partial\varphi_+/\partial x^n -\partial\varphi_T/\partial x^n}{\partial\varphi_+/\partial x^n +\partial\varphi_T/\partial x^n  } \quad \text{for $x^n=0$}.
\ee
As explained below \r{1.8}, we can get initial conditions for the subsequent transport equations, and then solve those transport equation. By \r{1.4}, the maps
\be{1.16a}
P_R: h\mapsto u_R|_{\R\times\Gamma}, \quad P_T: h\mapsto u_T|_{\R\times\Gamma}
\ee
are \PDO s of order $0$ with principal symbols equal to $b_R^{(0)}$, $b_T^{(0)}$ restricted to $\R\times\Gamma$, see \r{1.16}. We recall \r{1.23h} as well. 

We estimate next  the  amount of energy that is transmitted in  $\Omega_\text{ext}$. We will do it only based on the principal term in our parametrix. That corresponds to an estimate of the solution operator corresponding to transmission, up to compact operators, as we show below. 

A quick look at \r{1.16}, see also \r{1.5a} shows that $b_T^{(0)}>1$. This may look strange because we should have only a fraction of the energy transmitted, and the rest is reflected. There is no contradiction however because the energy is not proportional to the amplitude. % restricted to $\R\times\Gamma$. 

Let $u$ solve $(\partial_t^2-c^2\Delta)u=0$ in the bounded domain $U$ with smooth boundary for $t'\le t\le t''$ with some $t'<t''$. A direct calculation yields
\be{1.17}
E_U(\mathbf{u}(t''))=E_U(\mathbf{u}(t'))+2\Re \int_{[t',t'']\times \partial U}u_t\frac{\partial\bar  u}{\partial\nu}\,\d t\,\d S.
\ee
We will use this to estimate the energy of $u_T$ in   $\Omega_\text{ext}$.  Since the wave front set of $u_T$ is contained in some small neighborhood of the transmitted bicharacteristic,  we have smooth data for $t=0$. 
Therefore, if $t_2>t_1$ is fixed closed enough to $t_1$, we can apply \r{1.17} to a large ball minus $\Omega_\text{int}$ to get that modulo a compact operator applied to $h$,
\be{1.18}
E_{\Omega_\text{ext}}(\mathbf{u}_T(t_2)) \cong 2 \Re\int_{[0,t_2]\times \Gamma}\frac{\partial u_T}{\partial t}\frac{\partial \bar u_T}{\partial\nu}\,\d t\,\d S.
\ee
Therefore, 
\be{1.30}
E_{\Omega_{ext}}(\mathbf{u}_T(t_2)) \cong 2\Re (P_t u_T,N_{\text{ext,out}} u_T)=\Re (2P_T^*N_{\text{ext,out}}^*P_tP_Th,h),
\ee
where $(\cdot,\cdot)$ is the inner product in $\R\times\R^{n-1}$, and 
$P_t=\d/\d t$.

Apply similar arguments to $u_+$ in $\Omega_\text{int}$. Since the bicharacteristics leave $\Omega_\text{int}$, we have modulo smoother terms
\be{1.20}
0\cong E_{\Omega_\text{int}}(\mathbf{u}_+(0))+2 \Re\int_{[0,t_2]\times \Gamma}\frac{\partial u_+}{\partial t}\frac{\partial \bar u_+}{\partial\nu}\,\d t\,\d S.
\ee
Similarly we get, see again \r{1.16a},
\be{1.31}
E_{\Omega_\text{int}}(\mathbf{u}_+(0)) \cong -2\Re (P_t h,N_{\text{int,in}} h)= \Re (2N_{\text{int,in}}^*P_t h, h).
\ee
For the principal symbols of the operators in \r{1.30}, \r{1.31} we have 
\be{1.32}
\frac{\sigma_p(2P_T^*N_{\text{ext,out}}^*P_tP_T) }{ \sigma_p(2N_{\text{int,in}}^*P_t )}
= \frac{\partial\varphi_T/\partial\nu}{\partial\varphi_+/\partial\nu} \left(b_T^{(0)}\right)^2
 = \frac{4 (\partial\varphi_+/\partial\nu)( \partial\varphi_T/\partial\nu)}{(\partial\varphi_+/\partial\nu+\partial\varphi_T/\partial\nu)^2} .
\ee
Denote for a moment $a:= \partial\varphi_+/\partial\nu$, $b:= \partial\varphi_T/\partial\nu$. Then the quotient above equals $4ab/(a+b)^2\le 1$ that confirms that the reflected wave has less energy than the incident one. By \r{1.17_0}, $0\le b\le\gamma a$, $0<\gamma<1$. This easily implies
\be{1.23}
\frac{4ab}{(a+b)^2}\le \frac{4\gamma}{(1+\gamma)^2}<1.
\ee
Therefore, the expression in the middle represents an upper bound of the portion of the total energy that gets transmitted in the asymptotic regime when the frequency tends to infinity. To get a lower bound, assume in addition that $b\ge b_0>0$ and $a\le a_0$ for some $a_0$,  $b_0$,  i.e.,
\be{1.23m}
0<b_0< \frac{\partial\varphi_T}{\partial\nu}, \quad \frac{\partial\varphi_+}{\partial\nu}\le a_0.
\ee
Then
\be{1.23n}
\frac{4ab}{(a+b)^2}\ge \frac{4b_0^2/\gamma}{(1+\gamma)^2a_0^2}>0.
\ee
This is a lower bound  of the ratio of the high frequency energy that is transmitted. As we can  see, if the transmitted ray gets very close to a tangent one, that ratio tends to $0$. 

So far this is still not a proof of such a statement but just a heuristic argument. We will formulate and prove this below in Proposition~\ref{pr_energy}.

\subsection{Snell's Law}
Assume now that $(\xi^1)'$ is in the hyperbolic region for $c_\text{int}$ but not necessarily for $c_\text{ext}$. This corresponds to a ray hitting $\Gamma$ from the ``interior'' $\Omega_\text{int}$. There is no change in solving the eikonal equation for $\varphi_R$ but a real phase $\varphi_T$ does not exist if the expression under the square root in \r{1.27} is negative. This happens when $(\xi^1)'$ is in the elliptic region for $c_\text{ext}$.  
 Then there is no transmitted singularity in the parametrix. We analyze this case below. If $c_\text{int}>c_\text{ext}$, then  $(\xi^1)'$ that is in the hyperbolic region for $c_\text{int}$ by assumption, also falls into the hyperbolic region for the speed $c_\text{ext}$, i.e., there is always a transmitted ray. If $c_\text{int}<c_\text{ext}$, then existence of a transmitted wave depends on where $(\xi^1)'$ belongs w.r.t.\ $c_\text{ext}$.

Let $\alpha$ be the angle that $\xi^1 = \partial\varphi_+/\partial x^n$ makes with the (co)-normal represented by $\d x^n$, and let $\beta$ be the angle between  the latter and $\xi_T := \partial\varphi_T/\partial x^n$. We have
\be{1.17f}
|\xi'|=|\xi^1| \sin\alpha = c_\text{int}^{-1}\tau \sin\alpha, \quad 
|\xi'|=|\xi_T| \sin\beta = c_\text{ext}^{-1}\tau \sin\beta
\ee
By \r{1.17f},  we recover  Snell's law
\be{1.14}
\frac{\sin\alpha}{\sin\beta} = \frac{c_\text{int}}{c_\text{ext}},
\ee
Assume now that $c_\text{int}<c_\text{ext}$, which is the case where there might be no transmitted ray. 
Denote by 
\be{1.13c}
\alpha_0(x)= \arcsin(c_\text{int}(x)/c_\text{ext}(x))
\ee
the critical angle at any $x\in \Gamma$  that places $(\xi^1)'$ in the glancing manifold w.r.t.\ $c_\text{ext}$. Then the transmitted wave does not exist when $\alpha>\alpha_0$; more precisely we do not have a real phase function $\varphi_T$ in that case. It exists, when $\alpha<\alpha_0$. In the critical case $\alpha=\alpha_0$, this construction provides an outgoing ray tangent to $\Gamma$ that we are not going to analyze.

\subsection{The full internal reflection case} \label{sec_FIR}
Assume now that $(\xi^1)'$ is  in the elliptic region w.r.t.\ $c_\text{ext}$, then there is no transmitted singularity, but one can still construct a parametrix for the ``evanescent'' wave in  $\Omega_\text{ext}$; and there is a reflected ray. This is known as a full internal reflection. 
We give details below.

We proceed as above with one essential difference. There is no real valued solution $\varphi_T$ to the eikonal equation \r{1.24a} outside $\Omega_0$. Similarly to \r{1.27}, we get formally,

\be{1.50}
 \frac{\partial\varphi_T}{\partial\nu}  = 
\quad  \i\sqrt{ |\xi'|^2- c_\text{ext}^{-2}\tau^2}\quad 
\text{for $x^n=0$}. 
\ee
The choice of the sign of the square root is dictated by the requirement that the so parametrix \r{1.24} with $\sigma=T$ be exponentially decreasing away from $\Gamma$ instead of exponentially increasing. 

In general, the eikonal equation may not be solvable but one can still construct solutions modulo $O((x^n)^\infty)$. The same applies to the transport equations. One can show that the $O((x^n)^\infty)$ error does not change the properties of $u_T$ to be a parametrix.  
In particular, in \r{1.30} in this case one gets 
\be{1.51}
E_{\Omega_\text{ext}}(\mathbf{u}_T(t_2)) \cong 0,
\ee
because the principal term of $\partial\bar u_T/\partial\nu$ in \r{1.18} now is pure imaginary instead of being real. Moreover, $u_T$ is smooth in $\bar\Omega_\text{ext}$. 
Therefore, no energy, as far as the principal part only is considered, is transmitted to $\Omega_\text{ext}$. That does not mean that the solution vanishes there, of course.

\subsection{Glancing, gliding rays and other cases} We do not analyze the cases where $(\xi^1)'$ is in the glancing manifold w.r.t.\ to one of the speeds. We can do that because the analysis of those cases in not needed because of our assumptions guaranteeing no tangent rays. The analysis there is more delicate, and we refer to \cite{Taylor76, Petkov82,Petkov82II} for more details and examples. We do not analyze either the case where $(\xi^1)'$ is in the elliptic region with respect to either speed.  

\subsection{Justification of the parametrix} 
Denote by $\mathbf{u}_R=[u_R,\partial_t u_R]$, $\mathbf{u}_T=[u_T,\partial_t u_T]$ the approximate solutions constructed above, defined for $t$ in some neighborhood of $t_2$. Then $\mathbf{u}_R  = \mathbf{V}_Rh$, $\mathbf{u}_T  = \mathbf{V}_Th$, where $\mathbf{V}_{R,T}$ are the FIOs constructed above. Let $\mathbf{u}_+$ be the solution of \r{1'} defined above, with initial data $\mathbf{\Pi_+f}$ at $t=0$ having wave front set in a small neighborhood of $(x_0,\xi^0)$. The map $\Lambda_+: \mathbf{f}\mapsto u_+|_{\R\times\Gamma}=h$ is an FIO described in \cite{SU-thermo}. Then near $(t_1,x_1)$, 
\[
\mathbf{u}_R = \mathbf{V}_R\Lambda\mathbf{f}, \quad \mathbf{u}_T = \mathbf{V}_T\Lambda\mathbf{f},
\]
the former supported  in $\R\times\bar\Omega_\text{int}$, and the later in $\R\times \bar\Omega_\text{ext}$. So far we had two objects that we denoted by $u_+$: first, the parametrix of the solution of \r{1'} corresponding to the positive sound speed $c(x)|\xi|$; and the parametrix in $\R\times\bar\Omega_\text{int}$ for the incoming solution corresponding to boundary value $h$. When $h=\Lambda_+\mathbf{f}$, those two parametrices coincide up to a smooth term, as it is not hard to see (the second one is a back-projection and  is discussed in \cite{SU-thermo}, in fact). This justifies the same notation for them that we will keep.

Consider the parametrix $v_p := {u}_+ + {u}_R +  {u}_T$.  We can always assume that its support is in some small neighborhood of the geodesic that hits $\R\times\Gamma$ at $(t_1,x_1)$ and is tangent to $\xi^1$ there; and then reflects, and another branch refracts, see Figure~\ref{fig:skull1}. In particular, then $h$ has $t$-support near $t=t_1$, let us say that this included in the interval $[t_1-\eps,t_1+\eps]$ with some $\eps>0$. 
At $t=t_2:=t_1+2\eps$, let $x_2$ be the position of the reflected ray, and let $\xi^2$ be its unit co-direction. Then $\WF(u_R(t_2,\cdot))$ is in a small conic neighborhood of $(x_2,\xi^2)$. 

Let $\mathbf{v}(t,\cdot)  = e^{t\mathbf{P}}\mathbf{\Pi_+f}$ be the exact solution, with some fixed choice of the parametrix $Q^{-1}$ in the definition of $\mathbf{\Pi}_+$, properly supported.  Consider $w=v-v_p$ in $[0,t_2]\times\R^n$.  It satisfies

\begin{align}\label{1.24'}
(\partial^2_t-c^2\Delta)w|_{ [0,t_2]\times \bar\Omega_\text{int}}&\in C^\infty, 
&(\partial^2_t-c^2\Delta)w|_{ [0,t_2]\times \bar\Omega_\text{ext}}&\in C^\infty,\\
w|_{ [0,t_2]\times \Gamma_\text{ext}}-w|_{ [0,t_2]\times \Gamma_\text{int} } &\in C^\infty,   &\frac{\partial w}{\partial\nu}\big|_{ [0,t_2]\times \Gamma_\text{ext} }- \frac{\partial w}{\partial\nu}\big|_{ [0,t_2]\times \Gamma_\text{int}} &\in C^\infty. \label{1.25'}
\end{align}
On the other hand, for $0\le t\ll1$, $v$ is smooth. Let $\chi\in C^\infty(\R)$ be a function that vanishes in $(-\infty,\delta]$ and equals $1$ on $[2\delta,\infty)$, $0<\delta\ll1$. Then $\tilde w:= \chi(t)w(t,x)$ still satisfies \r{1.24'}, \r{1.25'} and also vanishes for $t\le0$. By \cite[Theorem~1.36]{Williams-transmission}, $\tilde w$ is smooth in $[0,t_2]\times \bar\Omega_\text{int}$, up to the boundary, and is also  smooth in $[0,t_2]\times \bar\Omega_\text{ext}$, up to the boundary. Therefore,
\be{1.26'}
\mathbf{v}(t,\cdot)= \mathbf{v}_p(t,\cdot) + \mathbf{K}_t\mathbf{f},
\ee
for any $t\in [0,t_2]$, where $\mathbf{K}_t$ is a compact operator in $\mathcal{H}$, depending smoothly on $t$. The operator $\mathbf{K}_t$ depends on $\mathbf{Q}$ as well. Therefore, the parametrix coincides with the exact solution up to a compact operator that is also smoothing in the sense described above.

\section{Proof of the main result}
We start with a lemma, that in principle is well known and is related to the classical Dirichlet principle.

\begin{lemma}\label{lemma_D}
Let $\Omega$ be a bounded set with smooth boundary. 
Define the map 
\[
\Pi_\Omega : H^1(\Omega) \longmapsto H_D(\Omega)
\]
as follows: $\Pi_\Omega u=u-\phi_u$, where $\phi_u$ is the solution to
\[
 \Delta \phi_u=0 \quad \text{in $\Omega$}, \quad \phi_u|_{\bo} = u|_{\bo}.
\]
Then 
\[
 \|\Pi_\Omega u\|_{H_D(\Omega)} \le \|u\|_{H_D(\Omega)},
\]
where $\|u\|_{H_D(\Omega)}$ is the Dirichlet norm of $u$ extended to functions in $H^1(\Omega)$ that may not vanish on $\bo$.

Set $\mathbf{\Pi}_{\Omega}\mathbf{u} = [\Pi_\Omega u_1,u_2]$. Then
\[
 E_\Omega(\mathbf{\Pi}_{\Omega}\mathbf{u})  \le E_\Omega(\mathbf{u}).
\]
\end{lemma}

\begin{proof}
In what follows, $(\cdot,\cdot)_{H_{D}(\Omega)}$ is the inner product in $H_{D}(\Omega)$, see \r{2.0H}, applied to functions that belong to $H^1(\Omega)$ but maybe not to $H_{D}(\Omega)$ because they may not vanish on $\bo$.  Since $u=\phi_u$ on $\bo$, and since $\phi_u$ is harmonic, we get
\[
(u-\phi_u,\phi_u)_{H_{D}(\Omega)}  = \int_\Omega \nabla (u-\phi_u)\cdot\nabla \bar\phi_u\,\d x
=0.
\]
Then
\[
\|u-\phi_u\|^2_{H_{D}(\Omega)} = \|u\|^2_{H_{D}(\Omega)} - \|\phi_u\|^2_{H_{D}(\Omega)}
\le \|u\|^2_{H_{D}(\Omega)}.
\]
The second statement of the lemma follows immediately from the first one.
\end{proof}

It is easy to see that $\Pi_\Omega$  is a projection, orthogonal with respect to  the product  
$(\cdot,\cdot)_{H_D(\Omega)}$. This products defines a seminorm  on $H^1(\Omega)$ only.

Before giving the proof of the main theorem, we will show that condition \r{2.1b} implies the following local energy decay estimate.

\begin{proposition}\label{pr_energy} 
Let \r{2.1b} be satisfied. Then there exists $0<\mu<1$ so that
\be{2.E}
E_\Omega(\mathbf{u}(T))\le \mu E_\Omega(\mathbf{f})
\ee
for any solution $\mathbf{u}$ with initial data $\mathbf{f}=[f_1,0]\in \mathcal{H}(\mathcal{K})$. 
\end{proposition}

\begin{proof}

We will use the geometric optics construction in section~\ref{sec_ans}. Let $\mathbf{X}=\text{diag}(X,X)$ be a zeroth order \PDO\ with small enough essential support  supported near  some $(x_0,\xi^0)\in T^* \mathcal{K} \setminus 0$. Let $\mathbf{v} := e^{t\mathbf{P}}\mathbf{\Pi_+ Xf}$ be the corresponding solution. Recall that $\mathbf{\Pi}_+$ restricts propagation of singularities to the positive wave speed $c(x)|\xi|$ only, while $\mathbf{X}$ localizes near $(x_0,\xi^0)$. Assume for now, that the ray through $(x_0,\xi^0)$ gives rise to both a reflected and a transmitted one, as in section~\ref{sec_r_t}.

We will use the energy computation in \r{1.30} and \r{1.31}, under the assumption \r{1.23m}. 
 According to those relations, \r{1.32} and \r{1.23n},
\be{2.3}
E_{ \Omega_\text{ext}}(\mathbf{u}_T(t_2))
-\mu' E_{\Omega_\text{int}}(\mathbf{u}_+(0))  = \Re (Mh,h),
\ee
with $M$ a \PDO\  of order $0$ with non-negative principal symbol, some $\mu'\in (0,1)$, and we use the notation in that section. By the G\.{a}rding inequality,
\be{2.3'a}
E_{\Omega_\text{ext}}(\mathbf{u}_T(t_2))\ge \mu' E_{\Omega_\text{int}}(\mathbf{u}_+(0))- C \|h\|^2_{H^{1/2}(\R\times\bo)} 
= \mu' E_{ \Omega_\text{int}}(\mathbf{\Pi_+ Xf}) 
- C \|h\|^2_{H^{1/2}(\R\times\bo)} .
\ee
The map $\Lambda_{1,+} : [f_1,0]\mapsto h$ is an FIO of order $0$  with a canonical relation of graph type discussed above. Similarly, the map $\Lambda_{2,+} : [0,f_2]\mapsto h$ is an FIO of order $-1$ with the same canonical relation. Then 
\[
\|h\|_{H^{1/2}(\R\times\bo)} \le C\| \mathbf{u}_+(0)\|_{H^{1/2}\oplus H^{-1/2}}.
\]
Therefore,
\be{2.3'}
E_{\Omega_\text{ext}}(\mathbf{u}_T(t_2))\ge \mu' E_{ \Omega_\text{int}}(\mathbf{\Pi_+ Xf}) 
- C\|  \mathbf{\Pi_+ Xf} \|^2_{H^{1/2}\oplus H^{-1/2}}.
\ee

Next, assume that the ray through $(x_0,\xi^0)$ gives rise to  a reflected ray only, as in  section~\ref{sec_FIR}. Then by \r{1.51},
\be{2.3d1}
E_{ \Omega_\text{ext}}(\mathbf{u}_T(t_2)) \le C\|  \mathbf{\Pi_+ Xf} \|^2_{H^{1/2}\oplus H^{-1/2}}.
\ee
By energy preservation, 
\be{2.3e}
E_{ \Omega_\text{ext}}(\mathbf{u}_R(t_2)) \ge    E_{ \Omega_\text{int}}(\mathbf{\Pi_+ Xf}) -   C\|  \mathbf{\Pi_+ Xf} \|^2_{H^{1/2}\oplus H^{-1/2}},
\ee
and we used the fact here that $\mathbf{u_+}(t_2,\cdot)$ can be obtain from the initial conditions by applying a smoothing operator.

We will apply those arguments now in a more general situation. Assume now that $(x_0,\xi^0)\in \mathcal{U}$, see \r{2.1a'}. By the definition of $\mathcal{U}$, there is a path starting from $(x_0,\xi^0)$ or from $(x_0,\xi^0)$  so that it consists of finitely many geodesic segments; and  at least one of them is a geodesic hitting $\Gamma$  transversely for the first time at some $x_1$ an angle strictly greater than the critical one $\pi/2-\alpha_0$ at $x_1$. Without loss of generality, we can assume that we have the former case: the path starts from  $(x_0,\xi^0)$. 
Then the next geodesic segment in that path is a transmitted one that is not tangent to $\Gamma$. The inequality \r{1.23m} then holds for the corresponding phase with some $b_0>0$. 

Assume that the essential support of the symbol $q$ is so small that the analysis in section~\ref{sec_FIR} applies near each full internal reflection. Then step by step, applying \r{2.3e} consecutively, we can reduce the analysis to the case where $x_1$ is the first point where the geodesic from $(x_0,\xi^0)$  hits $\Gamma$. Then we apply \r{2.3'} to estimate the amount of energy that has been transmitted on the other side of $\Gamma$. Let for a moment assume that  the transmitted ray  leaves $\Omega$ for time $t<T$ without further contact with $\Gamma$. Then for any $t>t_1$, not necessarily close to it, the map $\mathbf{u}_T(t_1,\cdot) \mapsto \mathbf{u}_T(t,\cdot) $ is an FIO  that represents the solution of the wave equation outside $\Omega_0$ corresponding to the positive wave speed. In particular, it is not affected by the presence of a jump at $\Gamma$.  
Therefore, an estimate equivalent to \r{2.3'} is preserved for such $t$ as well but now the energy of $\mathbf{u}(T)$ is concentrated outside $\bar\Omega$, up to smoothing terms. 
Let $\Omega_0$ be an open set with a smooth boundary so that $\mathcal{K}\subset \Omega_0\Subset\Omega\setminus\Gamma$. 
Then
\be{2.3.1n} 
\| \mathbf{X\Pi_+f}\|^2_{\mathcal{H}( \Omega_0)} \le 
C E_{\R^n\setminus \Omega} (\mathbf{u}(T))   +C\|  \mathbf{f} \|^2_{H^{1/2}\oplus H^{-1/2}}.
\ee
Above, we also modify $\mathbf{\Pi}_+$ and $\mathbf{X}$ by smoothing operators, if needed, so that $ \mathbf{X\Pi_+f}$  belongs to the indicated energy space, i.e., they vanish in $\R^n\setminus \Omega_0$. 

Consider now the case where the ray guaranteed by \r{2.1a'}, \r{2.1b} may hit $\Gamma$ again, even more than once, before leaving $\bar\Omega$ and not coming back. At each such event, there will be no loss of energy at high frequencies, as in \r{2.3d1}, or there will be positive portion of the transmitted high frequency  energy, as in \r{2.3'}. Then \r{2.3.1n} is still true. 

So far we did not use the assumption $f_2=0$. Now, since $\mathbf{f}=[f_1,0]$, we get $\mathbf{\Pi}_+\mathbf{Xf}=\frac12[f_1,\i QXf_1]$, see \r{42.4}.  By the ellipticity of $Q$, 
\be{2.3.1a} 
\|Xf_1\|_{H_D(\Omega_0)} \le 
C \|\mathbf{u}(T))\|_{H^1(\R^n\setminus \Omega) \oplus L^2( \R^n\setminus \Omega) ) }  +C\|f_1 \|_{H^{1/2}}.
\ee

By a compactness argument, in a conical neighborhood of $S^*\mathcal{K}$, we can take a finite  pseudo-differential partition of unity $1=\sum_j \chi_j$ 
of symbols of \PDO s $X_j$ localizing in conical neighborhoods of a finite number of points $(x_j,\xi^j)\in S^*\mathcal{K}$. Thus we get
\be{2.3.2} 
\|f_1\|_{H_D(\Omega_0)} \le 
C \|\mathbf{u}(T))\|_{H^1(\R^n\setminus \Omega) \oplus L^2( \R^n\setminus \Omega) ) }   +C\|f_1 \|_{H^{1/2}}.
\ee
Consider the bounded  map
\be{2.3.5}
H_D(\mathcal{K})\ni f_1 \longmapsto \mathbf{u}(T)\in {H^1(\R^n\setminus \Omega) \oplus L^2( \R^n\setminus \Omega) ) } .
\ee
We claim that it is injective. Indeed, assume that for some $\mathbf{f}\in \mathcal{K}$, for the corresponding $u$ we have
\[
u(T,x) = 0, \quad \mbox{for $x\not\in\Omega$}.
\]
By finite domain of dependence in $\R\times\R^n\setminus\Omega$, where $c=1$, we get 
\be{2.7}
u(t,x) = 0 \quad \mbox{when $\dist_{\rm e}(x,\Omega)>|T-t|$},%, \, 0\le t\le T$}.
\ee
where $\dist_{\rm e}$ stands for the Euclidean distance. 
One the other hand, we also have  
\be{2.8}
u(t,x) = 0 \quad \mbox{when $\dist_{\rm e}(x,\Omega)>|t|$}.%, \, 0\le t\le T$}.
\ee
Here, we applied finite domain of dependence argument outside $\Omega$, as well. Note that this is not sharp, because  at least when $c$ does not jump at $\Gamma$, then the Euclidean distance $\dist_{\rm e}(x,\Omega)$ can be replaced by the distance in the metric $c^{-2}dx^2$ between $\mathcal{K}$ and $x$. 

\nofigure{
\begin{figure}[t] % float placement: (h)ere, page (t)op, page (b)ottom, other (p)age
  \centering
  % file name: C:/Users/Plamen/Documents/My PCTeX Files/current/skull3.pdf
  \includegraphics[bb=0 0 508 457,width=3.07in,height=2.76in,keepaspectratio]{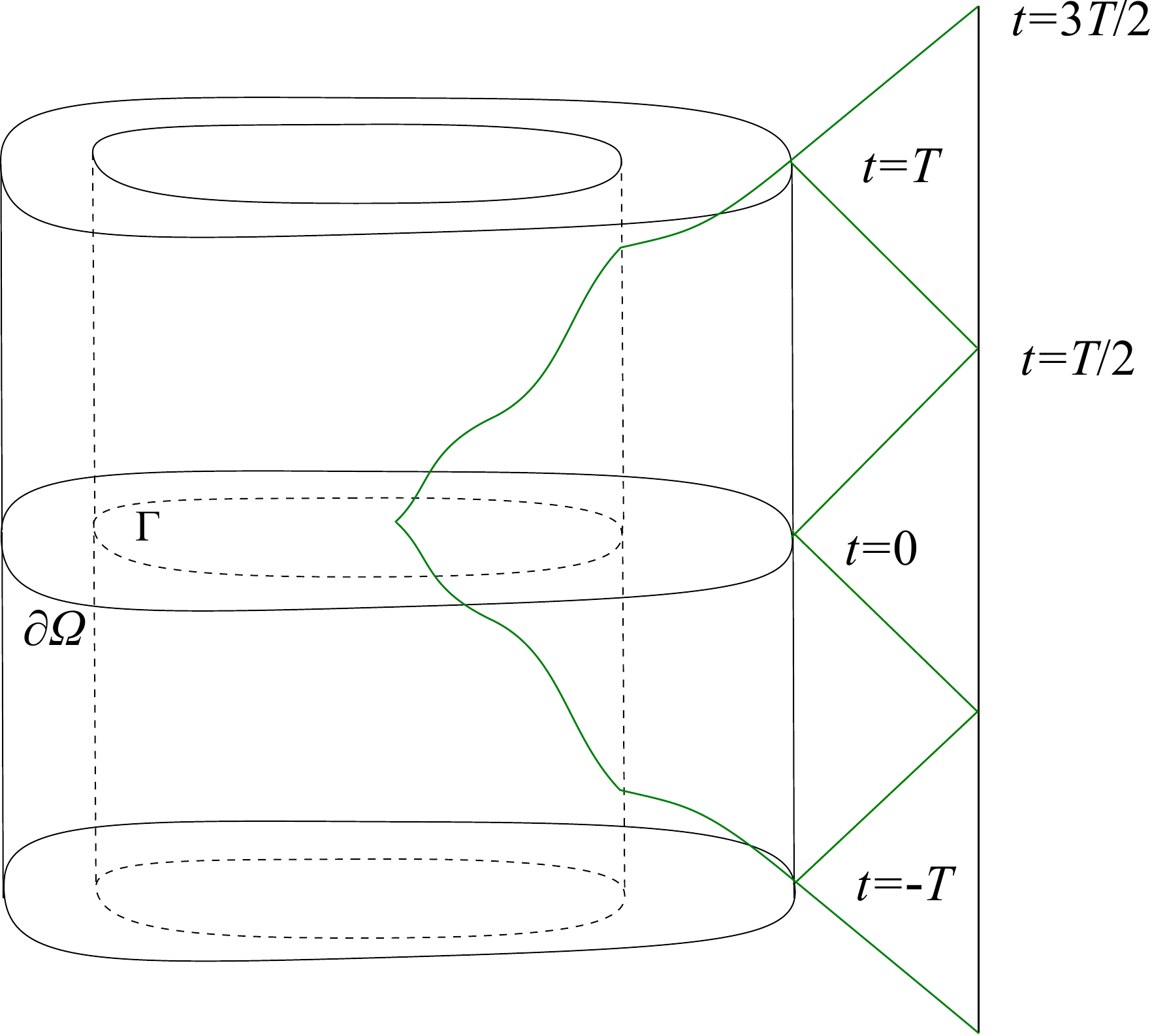}
 \caption{}
  \label{fig:skull3.pdf}
\end{figure}
}

Therefore,
\be{2.8a}
u(t,x) = 0 \quad \mbox{when $\dist_{\rm e}(x,\Omega)>T/2, \;  -T/2\le t\le 3T/2$}.
\ee
The solution  $u$ extends to an even function of $t$ that is still a solution of the wave equation because $f_2=0$. Then one gets that \r{2.8a} actually holds for $|t|<3T/2$. 
We will conclude next  by the unique continuation   Theorem~\ref{thm_uq} that $u=0$ on $[0,T]\times \Omega$, therefore, $f_1=0$.  

To this end, notice fist that from John's theorem (equivalent to Tataru's unique continuation result \cite[Theorem~2]{SU-thermo} in the Euclidean setting), we get  $u=0$ on $[-T,T]\times \R^n\setminus \Omega$. Fix $x_0\in\Omega$. Then there is a piecewise smooth curve starting at $x_0$ in direction either $\xi^0$ or $-\xi^0$, where  $\xi^0$  is arbitrary and fixed, of length less than $T$ that reaches $\bo$. This means that $\dist(x_0,\R^n\setminus\Omega)<T$, see Definition~\ref{def_dist}. Then by Theorem~\ref{thm_uq}, $u(0,\cdot)=0$ near $x_0$. Since $x_0$ was arbitrary, we get $f=0$. This completes the proof of the claim that \r{2.3.5} is injective. 

Now, by \cite[Proposition~5.3.1]{Taylor-book1}, since the inclusion $H_D(\mathcal{K}) \hookrightarrow H^{1/2}$ is compact, 
there is an estimate as in \r{2.3.2}, with a different $C$, with the last term missing, i.e.,
\be{2.3.6} 
\|f_1\|_{H_D(\mathcal{K})} \le C\| \mathbf{u}(T) \|_{H^1(\R^n\setminus \Omega)\oplus L^2(\R^n\setminus \Omega)} .
\ee
By finite speed of propagation, $\mathbf{u}(T)$ is supported in some large ball $B$. 
Apply the Poincar\'e inequality \cite{Gunther-book} to conclude that 
\be{2.3.6P} 
\|f_1\|^2_{H_D(\mathcal{K})} \le CE_{\R^n\setminus \Omega}(\mathbf{u}(T)) .
\ee

Then, using Lemma~\ref{lemma_D}, we get
\[
\begin{split}
E_\Omega(\mathbf{u}(T))& = E_{\R^n}(\mathbf{u}(T)) - E_{\R^n\setminus \Omega}(\mathbf{u}(T))
 \\
  &\le E_{\mathcal{K}}(\mathbf{f}) - E_{\mathcal{K}}(\mathbf{f})/C\\
   &= (1-1/C^2)E_{\mathcal{K}}(\mathbf{f}),
\end{split}
\]
and this completes the proof.
\end{proof}

\begin{proof}[Proof of Theorem~\ref{thm2.1}]
Let $\mathbf{f}\in \mathcal{H}(\Omega)$ first. Let $w$ solve
\begin{equation}   \label{2.3d}
\left\{
\begin{array}{rcll}
(\partial_t^2 -c^2\Delta)w &=&0 &  \mbox{in $(0,T)\times \Omega$},\\
w|_{[0,T]\times\partial\Omega}&= &0,\\
w|_{t=T} &=& u|_{t=T}-\phi,\\ \quad w_t|_{t=T}& =&u_t|_{t=T},\\
\end{array}
\right.               
\end{equation}
where $u$ solves \r{1'} with a given $\mathbf{f}\in \mathcal{H}$, and $\phi$ is as in \r{l2} with $h=\Lambda f$. We can also write
\be{2.3.7}
[w,w_t]|_{t=T} = \mathbf{\Pi}_{\Omega}\left( [u,u_t]|_{\{T\}\times\Omega}\right),
\ee
where $\mathbf{\Pi}_{\Omega} : H^1(\Omega) \to\mathcal{H}(\Omega)$ is the  projection introduced in Lemma~\ref{lemma_D}. Then
\be{2.3.7d}
E_\Omega(\mathbf{w}(T)) 
 \le E_{\Omega}(\mathbf{u}(T)).
\ee

Let $v$ be the solution of \r{l2} with $h=\Lambda f$. Then $v+w$ solves the same initial boundary value problem in $[0,T]\times\Omega$ that $u$ does (with initial conditions at $t=T$), therefore $u=v+w$. Let now $\mathbf{f}=[f_1,0]\in \mathcal{H}(\mathcal{K})$. 
 Restrict ${u}$  to $t=0$ and  project to ${H}_D(\mathcal{K})$ to get
\[
f_1= \Pi_{\mathcal{K}}A_1\Lambda_1  f_1 +{\Pi}_{\mathcal{K}} {w}(0,\cdot).
\]
Set 
\[
{K} {f}_1 = {\Pi}_{\mathcal{K}}{w}(0,\cdot) .
\]
It remains to show now that ${K}$ extends to a contraction. 

By \r{2.3.7d} and Proposition~\ref{pr_energy}, the energy  of the initial conditions in \r{2.3d} satisfies the inequality
\be{2.4}
E_\Omega(\mathbf{w}(T)) \le \mu \|f_1\|^2_{H_D}.
\ee
Since the Dirichlet boundary condition is energy preserving, we get 
\be{2.4a}
E_{\Omega}(\mathbf{w}(0)) =    E_{\Omega}(\mathbf{w}(T))\le \mu \|f_1\|^2_{H_D}. 
\ee
Therefore, 
\be{2.5}
\|Kf_1\|^2_{H_{D}(\mathcal{K})} \le  E_{\mathcal{K}}(\mathbf{\Pi}_{\mathcal{K}}{w}(0))\le 
E_{\mathcal{K}}(\mathbf{w}(0))\le  E_{\Omega}(\mathbf{w}(0))\le \mu \|f_1\|^2_{H_D}.
\ee
This completes the proof. 
\end{proof}

\section{Unique continuation}\label{sec_FS}
\begin{definition}\label{def_dist}
Given $x$, $y$ in $\R^n$, let $\dist(x,y)$ be the infimum of the length of all %piecewise 
smooth curves connecting $x$ and $y$, intersecting $\Gamma$ transversally at each common point, where the length is measured in the metric $c^{-2}dx^2$.
\end{definition}

\medskip
The next result  is a   unique continuation theorem. The proof is based on the smooth case.

\begin{theorem}\label{thm_uq}
Let $x_0\not\in \bar\Omega$ and $T>0$. 
Assume that    $[u,u_t]\in C(\R;\, \mathcal{H}_\text{\rm loc})$ 
and $u$ satisfies 
\[
(\partial_t^2-c^2\Delta)u=0,
\]
in a neighborhood of the set 
\be{FS1a}
U = \{(t,x);  \; |t|+\dist(x_0,x) \le T\}
\ee 
and vanishes in a neighborhood of $[-T,T]\times \{x_0\}$.
Then 
\be{FS1}
u(t,x)=0 \quad \mbox{in $U$}.
\ee
\end{theorem}

\begin{proof}
When $c$ is smooth, the theorem follows from Tataru's unique continuation result \cite{tataru95, tataru99}, see \cite{SU-thermo}. In the case we consider, we will base our proof on the smooth case. 

We show next that one can assume without loss of generality that $[u,u_t]\in C(\R;\,D(\mathbf{P})) $, where $D(\mathbf{P})$ stands for the domain of $\mathbf{P}$ equipped  with the graph topology.  
Assume first that $u$ solves the wave equation in the whole $\R\times\R^n$.   
By finite speed of propagation \cite{Williams-transmission}, we can always assume that $u$ has compactly supported Cauchy data for $t=0$. Then 
$\mathbf{u}=e^{t\mathbf{P}}\mathbf{f}$, with some $\mathbf{f}\in\mathcal{H}$. 
Take $\phi\in C_0^\infty(\R)$ and consider the convolution $u_\phi=u*_t\phi$ of $u$ with $\phi$ in the $t$ variable. Then $u_\phi$ is a solution of the wave equation with support close enough to that of $u$, if $\supp\phi$ is small enough; so if we prove the theorem for $u_\phi$, that would be enough. The initial value for $\mathbf{u_\phi}=[   u_\phi,\partial_t u_\phi]$ at $t=0$ is given by $\mathbf{u}_\phi(0,\cdot) = \int \mathbf{u}_\phi(-s,\cdot)\phi(s)\,\d s$. Apply the generator $\mathbf{P}$ to that to get 
\[
\mathbf{P}\mathbf{u}_\phi(0,\cdot) = \int \mathbf{P}\mathbf{u}_\phi(-s,\cdot)\phi(s)\,\d s= 
 \int \mathbf{u}_\phi(-s,\cdot)\phi'(s)\,\d s\in \mathcal{H}.
\]
Therefore, the initial condition now is in the domain of $\mathbf{P}$, 
and this proves our claim in this case.  Note that we can put any power of $\mathbf{P}$ there and get the same conclusion. 

Assume now that $u$ solves the wave equation  near $U$ only, as in the theorem. Let $u_\phi$ be as above. Fix $(t_1,x_1)\in U$ and let $V := \{(t,x);\; |t-t_1|<\eps, \, |x-x_1|<\eps\}$ with some $0<\eps\ll1$.
 By the non-sharp finite speed of propagation result in  \cite{Williams-transmission}, $u|_V$ depends on $u(t_0,\cdot)$, $t_0 := t_1-2\eps$,  restricted to a neighborhood of $x_1$ of size $O(\eps)$, assuming that all characteristic cones for the sound speed $c_0 :=\max c$  with vertices in $V$ and above $t=t_0$, lie entirely in $U$. The latter is true of $\eps\ll1$. Therefore, we can cut $u(t_0,\cdot)$ appropriately to make the support compact and contained in $U\cap \{t=t_0\}$, and then use it as initial data at $t=t_1-2\eps$. Then we write $\mathbf{v}(t,\cdot)= e^{\i (t-t_0)\mathbf{P}}\mathbf{g}$, where $\mathbf{g}$ is the localized Cauchy data at $T=t_0$ as above. If the cutoff is chosen appropriately, then $\mathbf{u}=\mathbf{v}$ on $V$ and we can apply the argument of the preceding paragraph.

We continue  with the observation that if the set $\{|t-t_0|+\dist(x_0,x)\le T_0\}$ for some $t_0$, $T_0$, does not intersect $\R\times \Omega_0$, then one can apply the ``smooth version'' of the theorem to conclude that any solution in a neighborhood of that set that vanishes near $[t_0-T_0,t_0+T_0]\times \{x_0\}$ also vanishes in that set. 

Next, note that it is enough to show that $u=0$ at one of the ``extreme'' points of $U$. More precisely, let $y$ be such that $\dist(x_0,y)=T$. Assume that we can prove that $u=0$ near $(t,x) = (0,y)$. Let $(\hat t, \hat x)$ be any point in $U$, and consider the set
\be{FS2}
|t-\hat t|+\dist(x_0,x)\le \hat T := \dist(x_0,\hat x).
\ee
It is included in $U$, and if we know to prove the result for the ``extreme'' points, then we would get that $u=0$ near $(\hat t, \hat x)$.

We divide the rest of the proof into several steps. Let $\Omega_0$ be interior of the complement of the maximal  unbounded connected component of $\R^n\setminus\Gamma$. Let $\Gamma_0=\bo_0\subset\Gamma$. Given $x$, $y$ both    
in $\R^n\setminus \Omega_0$, 
let $\dist_0 (x,y)$ be the infimum of the length of all %piecewise 
smooth curves in  $\R^n\setminus \Omega_0$, connecting $x$ and $y$, intersecting $\Gamma$ transversally at each common point, where the length is measured in the metric $c^{-2}dx^2$.

Choose and fix $y\in\Gamma_0$ so that
\be{FS3}
\ell := \dist_0(y,x_0)< T. 
\ee
If there is no such $y$, then the theorem follows by its ``smooth version''.  
Note that $\dist(y,x_0)\le \dist_0(y,x_0)$, therefore
\be{FS3a}
[-T+\ell,T-\ell]\times \{y\}\Subset  U.
\ee
Fix $0<\delta\ll1$.  Let $ y_\delta$ be the point in $\R^n\setminus \Omega_0$ that lies on the normal geodesic through $y$, and is at distance $\delta$ from $y$. Let $[0,1]\ni t\mapsto \gamma(t)$ be a smooth curve in $\R^n\setminus \bar\Omega_0$ connecting $x_0=\gamma(0)$ and $y_\delta=\gamma(1)$ of length (in the metric $c^{-2}\d x^2$)  $\ell_\delta$ not exceeding $\ell+3\delta$. It is easy to see that this can be done. Assume also that $\delta$ is so small that $\ell_\delta<T$, see \r{FS3}. Let $0<\eps=\dist_0(\gamma,\Gamma)/2$.

Let $x_1$, $x_2,\dots x_k$ be finitely many points on $\gamma$, corresponding to increasing values of $t$, so that each segment of $c(t)$ with endpoints $(x_0,x_1)$, $(x_1,x_2)$, \dots, $(x_k,y_\delta)$ is of length, that we denote by $\ell_{x_0x_1}$, etc., at most $\eps$. Apply the ``smooth'' version of the  theorem for the set 
\[
U_{\hat t, x_1} := \left\{(t,x);\;|t-\hat t|+\dist_0(x,x_0) \le \ell_{x_0x_1}\right\}
\]
with $\hat t$ such that the interval $\{t;\; |t-\hat t|\le \ell_{x_0x_1}\}$ is included in $[-T,T]$. Such $\hat t$ clearly exist and are described by $|\hat t|\le T-\ell_{x_0x_1}$. Next, for any $\hat t$ as above, we have $U_{\hat t, x_1}\subset U\cap (\R\times \R^n\setminus \bar\Omega)$.  Therefore, we can apply the ``smooth'' version of the theorem to conclude that $u=0$ in a neighborhood of $(\hat t, x_1)$ for any such $\hat t$. In other words, $u=0$ in a neighborhood of $[-T+\ell_{x_0x_1}, T-\ell_{x_0x_1}]\times \{x_1\}$. We can repeat those arguments to get  $u=0$ in a neighborhood of $[-T+\ell_{x_0x_1}+\ell_{x_1x_2}, T-\ell_{x_1x_2}-\ell_{x_1x_2}]\times \{x_2\}$, etc. Thus we get
\[
u=0\quad\text{near}\quad [-T+\ell_\delta, T-\ell_\delta]\times \{y_\delta\}.
\]
Take the limit $\delta\to0$ now to conclude that for any $\ell'>\ell$, $u=0$ on an open set containing $[-T+\ell', T-\ell']\times\{x;\;0<\dist_0(x,y)\ll1\}$. Perturb $y$, and use the fact that our assumptions on $T$ are open, i.e., we can perturb $T$ a bit as well, to get
\be{FS4}
u=0\quad \text{in an one-sided neighborhood of $[-T+\ell, T-\ell]\times \{y\}$ in $\R\times\R^n\setminus\Omega_0$},
\ee
and is a solution there. 
By the assumed regularity of $u$,  taking the trace of $u$ and its normal derivative on $\Gamma$, from outside (i.e., in $\R^n\setminus \Omega_0$)   is well defined; and this trace is zero near $ [-T+\ell, T-\ell]\times \{y\}$. By the transmission conditions \r{1.3}, the interior traces of $u$ and its normal derivative vanish there as well. Let now $\tilde c$ be any smooth extension of $c|_{\Omega_0}$ to $\R^n$. Then 
\be{FS5}
(\partial_t^2-\tilde c^2(x)\Delta )u=0\quad \text{in a two-sided neighborhood of $[-T+\ell, T-\ell]\times \{y\}$ in $\R\times\R^n$}.
\ee
By a two-sided neighborhood, we mean a normal one (an open set) --- we only used that term to emphasize the difference with the neighborhood in \r{FS4}. Since $u$ vanishes in a neighborhood of $[-T+\ell, T-\ell]\times \{y\}$ in the exterior, we apply the ``smooth'' version of the theorem for the speed $\tilde c$ near sets of the kind $[\hat t-\delta,\hat t+\delta]\times \{y\}$, $0<\delta\ll1$, $|\hat t|<T-\ell-\delta$ to conclude that $u$ also  vanishes in the neighborhood in \r{FS5}, after we shrink it if needed by $O(\delta)$. We can use our freedom to vary $T$ a bit to conclude that
\be{FS6}
u=0\quad \text{in a two-sided neighborhood of $[-T+\ell, T-\ell]\times \{y\}$ in $\R\times\R^n$}.
\ee
If small enough, that neighborhood is in $U$, see \r{FS3a}. Notice that we actually proved this property for any $y\in \R^n\setminus\Omega_0$ with $\ell$ as in \r{FS3}, not only for $y\in \Gamma$ but the latter case only requires the use of the transmission condition. 

For a fixed $a\in\Omega_0$, close enough to $\Gamma_0$, let $\Omega_1$ be the maximal open connected component of $\R^n\setminus\Gamma$ containing $a$. In $\Omega_1$, we define a distance function $\dist_1$ by minimizing over all smooth curves, transversal to $\Gamma$, that stay in $\bar\Omega_1$. 
Let $T_1=T-\ell>0$, see \r{FS3}, and choose $y_1\in \Gamma$, if possible,  so that 
\be{FS7}
\ell_1 := \dist_1(y_1,y)< T_1.
\ee
Use the arguments above with $T$ replaced by $T_1$, $x_0$ replaced by $y$ before, and $y$ replacing $y_1$. The only difference is that we work in $\Omega_1$ now. We then get
\be{FS7'}
u=0\quad \text{in a two-sided neighborhood of $[-T+\ell+\ell_1, T-\ell-\ell_1]\times \{y_1\}$ in $\R\times\R^n$}.
\ee
Let now $z$ be a point with $\dist(z,x_0)< T$. Choose $\delta>0$ so that there is a smooth curve $\gamma$ of length $\dist(z,x_0)+\delta$  crossing $\Gamma$ transversally each time, and connecting $x_0$ and $z$. Let also $\delta\ll1$ so that $\dist(z,x_0)+\delta< T$. 
By a compactness argument, $\gamma$ will cross $\Gamma$ finitely many times. Apply the argument above for each segment either in $\bar\Omega_0$ or in $\R^n\setminus\Omega_0$ to get 
\[
u=0 \quad\text{near $ [-T-\dist(z,x_0)-\delta, T+\dist(z,x_0)+\delta]\times \{z\} $}.
\]
This shows that $u=0$ in the interior of $U$. Since our assumptions allow us to increase $T$ slightly, $u=0$ near $U$ as claimed.
\end{proof}

%========================

%\newpage
\section{We can get the whole Cauchy data on $[0,T]\times\bo$}
We will show here that knowing $\Lambda \mathbf{f}$, one can recover the Neumann derivative of the solution at $[0,T]\times\bo$ as well. This is done by applying a non-local \PDO\ to $\Lambda \mathbf{f}$. This is known in principle, and we do not use it in our proofs. It reveals links to Control Theory however. 

We will define first the outgoing DN map. Given $g\in C_0^\infty([0,\infty)\times \bo)$, let $w$ solve the exterior mixed problem with $c=1$:
\begin{equation}   \label{1w}
\left\{
\begin{array}{rcll}
(\partial_t^2 -\Delta)w &=&0 &  \mbox{in $(0,T)\times \R^n$},\\
w|_{[0,T]\times\bo}&=&g,\\
w|_{t=0} &=& 0,\\ \quad \partial_t w|_{t=0}& =&0, 
\end{array}
\right.               
\end{equation}
Then we set
\[
Ng = \frac{\partial w}{\partial\nu} \Big|_{[0,T]\times\bo}.
\]
By  \cite{LasieckaLT}, for $g\in H^1_{(0)}([0,T]\times \bo)$, we have $[w,w_t]\in C([0,T);\;  \mathcal{H})$; therefore, 
\[
N: H^1_{(0)}([0,T]\times \bo) \to C([0,T]\times H^\frac12(\bo))
\]
is continuous. Note  that the results in  \cite{LasieckaLT} require the domain to be bounded but by finite domain of dependence we can remove that restriction in our case.  We also refer to \cite[Proposition~2]{finchPR} for  a sharp domain of dependence result for exterior problems.

\begin{lemma}\label{lemma_C}
Let $u$ solve \r{1'} with $\mathbf{f}\in \mathcal{H}$ compactly supported in $\Omega$. Assume that $c=1$ 
outside $\Omega$. Then for any $T>0$,  $\Lambda\mathbf{f}$ determines uniquely $u$ in $[0,T]\times \R^n\setminus\Omega$ and the normal derivative of $u$ on $[0,T] \times \bo$ as follows:

(a) The solution  $u$ in $[0,T]\times \R^n\setminus \Omega$ coincides with the solution of \r{1w} with $g=\Lambda\mathbf{f}$,

(b) We have 
\be{N1}
\frac{\partial w}{\partial\nu}\Big|_{[0,T]\times\bo} = N\Lambda f.
\ee
\end{lemma}

\begin{proof}
Let $w$ be the solution of \r{1w} with $g=\Lambda\mathbf{f}$. The latter is in $H^1_{(0)}([0,T]\times \bo)$, see the paragraph after \r{l1}. Let $u$ be the solution of \r{1'}. Then $u-w$ solves the unit speed wave equation in $[0,T]\times \R^n\setminus\Omega$ with zero Dirichlet data  and zero initial data. Therefore, $u=w$ in $[0,T]\times \R^n\setminus\Omega$. 
\end{proof}

\begin{remark}
Note that $c=1$ outside $\Omega$ was not a necessary assumption. 
\end{remark}

\bibliographystyle{abbrv}
%\bibliography{myreferences,TAT}
%

\end{document}